\documentclass{gOMS2e}

\usepackage{epstopdf}
\usepackage{subfigure}
\usepackage{amsfonts}
\usepackage{amsthm}
\usepackage{color}
\usepackage{amssymb}
\usepackage{amsmath}
\usepackage{geometry}
\usepackage{fancyhdr}
\usepackage{graphicx}
\usepackage{float}
\usepackage{multirow}
\usepackage[table,xcdraw]{xcolor}

\newcommand{\bit}{\begin{itemize}}
\newcommand{\eit}{\end{itemize}}
\newcommand{\be}{\begin{equation}}
\newcommand{\ee}{\end{equation}}
\newcommand{\Rm}{\mathbb{R}}

\theoremstyle{plain}
\newtheorem{proposition}{Proposition}[section]
\newtheorem{theorem}[proposition]{Theorem}

\newtheorem{algorithm}[proposition]{Algorithm}

\theoremstyle{definition}
\newtheorem{definition}[proposition]{Definition}

\theoremstyle{remark}
\newtheorem{remark}[proposition]{Remark}

\begin{document}



\title{An Efficient Augmented Lagrangian Method for Support Vector Machine}

\author{
\name{Yinqiao Yan$^{\ast}$\thanks{$^\ast$School of Mathematics and Statistics, Beijing Institute of Technology, Beijing, 100081, P. R. China}
and Qingna Li$^\dagger$\thanks{$^\dagger$Corresponding author. This author's research is supported by NSFC 11671036. School of Mathematics and Statistics/Beijing Key Laboratory on MCAACI, Beijing Institute of Technology, Beijing, 100081, P. R. China. Email: qnl@bit.edu.cn}}
}

\maketitle

\begin{abstract}
{\bf Abstract.} Support vector machine (SVM) has proved to be a successful approach for machine learning. Two typical SVM models are the L1-loss model for support vector classification (SVC) and $\epsilon$-L1-loss model for support vector regression (SVR). 
Due to the nonsmoothness of the  L1-loss function in the two models, most of the traditional approaches focus on solving the dual problem. In this paper, we propose an augmented Lagrangian method for the L1-loss model, which is designed to solve the primal problem. By tackling the non-smooth term in the model with Moreau-Yosida regularization and the proximal operator, the subproblem in {augmented Lagrangian method} reduces to a non-smooth linear system, which can be solved via the  quadratically convergent {semismooth Newton's method}. Moreover, the high computational cost  in semismooth Newton's method can be significantly reduced by exploring the sparse structure in the generalized Jacobian. Numerical results on various datasets in LIBLINEAR show that the proposed method is competitive with the most popular solvers in both speed and accuracy.
\end{abstract}


%
%
%
%
%
%
%

\section{Introduction}
Support vector machine (SVM) has proved to be a successful approach for machine learning.  {Support vector classification (SVC) and support vector regression (SVR) are two main types of support vector machines (SVMs). Support vector classification is a classic and well-performed learning method for two-group classification problems \cite{yj9}, whereas support vector regression is a learning machine extended from SVC by Boser {et al.} \cite{yj4}.
	{Given a training dataset, the  learning algorithms for SVC} can be used to find a maximum-margin hyperplane that divides all the training examples into two categories. It obtains the prediction function based on only a subset of support vectors. For SVR, instead of minimizing the training error, support vector regression minimizes the generalization error bound \cite{SVRreview} with a maximum-tolerance $\epsilon$, below which we do not have to compute the loss.}  {Two typical models} are the L1-loss SVC model and the $\epsilon$-L1-loss SVR model.

For {L1-loss SVC}, due to the nonsmoothness of the hinge loss function, traditional ways to deal with the hinge loss function is to introduce {slack} variables to formulate the problem as an optimization problem with a  smooth function and linear inequality constraints. Most approaches are proposed to solve the dual problem. For example, Platt \cite{Platt1998SMO} presented {a} sequential minimal optimization (SMO) algorithm. It deals with the dual problem by breaking the large quadratic programming(QP) optimization problem  into a series of small quadratic programming problems. The fast APG (FAPG) method \cite{fapg} is another widely used method to solve the QP problem with linear and bounds constraints.
Joachims proposed  ${\rm SVM^{light}}$ \cite{Joachims99a} and ${\rm SVM^{perf}}$ \cite{Joachims2006} respectively. ${\rm SVM^{light}}$ {is based on a generalized version of the decomposition strategy. In each iteration, all the variables of the dual problem are divided into two sets. One is the set of free variables and the other is the set of fixed variables.} ${\rm SVM^{perf}}$ uses a cutting-plane algorithm for training structural classification SVMs and ordinal regression SVMs.
Smola {et al.} \cite{Smola2008} improved ${\rm SVM^{perf}}$ by applying bundle methods and {it} performed well for large-scale data sets.
An exponentiated gradient (EG) method was proposed by Collins {et al.} \cite{Collins2008}. The EG method is based on exponentiated gradient updates and can be used to solve both the log-linear and max-margin optimization problems. 
Hsieh {et al.} \cite{LinLiblinear} proposed a dual coordinate descent (DCD) method for linear SVM to deal with large-scale sparse data. Methods aiming to  solve the primal form of the problem include the stochastic gradient descent method (SGD)   and {its different} variants
{such as averaged SGD (ASGD) \cite{Zhang2004}}. Shalev-Shwartz {et al.} \cite{Shwartz2007} proposed the Pegasos algorithm by introducing subgradient of the approximation for the objective function to cope with the non-differentiability of the hinge loss {function}. {They} considered a different procedure for setting the step size and included gradient-projection approach as an optional step. NORMA method proposed by Kivinen {et al.} \cite{NORMA} is also a variant of stochastic subgradient method based on the kernel expansion of the function. 
Inspired by Huber loss, Chapelle \cite{ChapelleHuber} used a differentiable approximation of L1-loss function. Tak${\rm \acute{a}\check{c}}$ {et al.} \cite{minibatchSGD} introduced the mini-batch technique in Pegasos algorithm to guarantee the parallelization speedups.
Recently, Chauhan et al. \cite{review18} presented a review on linear SVM and concluded that SVM-ALM proposed by Nie \cite{nie14} was the fastest algorithm which was applied to the L$_p$-loss primal problem by introducing the augmented Lagrangian method, and LIBLINEAR \cite{liblinear} is the most widely used solver which applied DCD to solving L2-regularized dual SVM. Niu et al. \cite{peipei2019} proposed SSsNAL method to solve the large-scale SVMs with big sample size by using the augmented Lagrangian method to deal with the dual problem.

For {the} $\epsilon$-L1-loss SVR model, similar to SVC, a new dual coordinate descent (DCD) method was proposed by Ho and Lin \cite{DCDSVR} for linear SVR. Burges {et al.} \cite{quadProSVR} applied an active set method to solving the dual quadratic programming problem. Smola \cite{smola96SVR} introduced a primal-dual path method to solve the dual problem. Similarly, stochastic gradient descent methods have good performance in solving large-scale support vector regression problem.

On the other hand, in optimization community, there are some recent progress on methodologies and techniques to deal with non-smooth problems. A typical tool is the semismooth Newton method, which is to solve non-smooth linear equations \cite{QiSun93}. Semismooth Newton's method has been successfully applied in solving various model optimization problems including nearest correlation matrix problem \cite{QiSun06}, nearest Euclidean distance matrix problem {\cite{Qi2013A}},  and so on {\cite{complementProb,LiQi2012,Qi2014A}. Moreover, Zhong and Fukushima \cite{multiSVM} use semismooth Newton's method to solve the multi-class support vector machines.} Recently, it is used to solve L2-loss SVC and $\epsilon$-L2-loss SVR \cite{YinLi2019}. For optimization problems including non-smooth terms in objective functions, Sun and his collaborators proposed different approaches based on the famous  Moreau-Yosida regularization. For example, to deal with the well-known LASSO problems, a highly efficient semismooth Newton augmented Lagrangian method is proposed in \cite{LiSun2017}. Similar technique is used to solve  the OSCAR and SLOPE {models},  as well as {convex} clustering \cite{Sun2017,Sun2018}.  To deal with two non-smooth terms in objective functions, an ABCD (accelerated block coordinate descent) framework \cite{ABCD} is proposed with the symmetric Gauss-Seidel technique embedded. The ABCD approach was applied to solve the Euclidean distance matrix model for protein molecular conformation {in} {\cite{ZhaiLi2019}}. 
In fact, the augmented Lagrangian method is quite popular and powerful to solve constraint optimization problems. With semismooth Newton's method as a subsolver, it is able to deal with various problems with non-smooth terms. The famous SDPNAL+ \cite{Zhao2009,SDPNAL,SDPNAL+,SDPNAL+matlab} is designed under the framework of augmented Lagrangian method.

Based on the above observations, a natural question arises. Given the fact that semismooth Newton's method has been used to solve L2-loss SVC and SVR, is it possible to solve the corresponding L1-loss models by making use of the modern optimization technique and approaches to tackle the non-smooth term? It is this question that motivates the work in our paper.

The contribution of the paper is as follows. Firstly, we propose an {augmented Lagrangian method} to solve the primal form of the L1-loss model for SVC and $\epsilon$-L1-loss model for SVR. The challenge of the nonsmoothness is tackled with  Moreau-Yosida regularization. Secondly, we apply semismooth Newton's method to solve the resulting subproblem. The quadratic convergence rate for semismooth Newton's method is guaranteed. Moreover, by exploring the sparse structure of the generalized Jacobian, the high computational complexity for semismooth Newton's method can be significantly reduced. Finally, extensive numerical tests and  datasets in LIBLINEAR demonstrate the fast speed and impressive accuracy of the method.

The organization of the paper is as follows. In Section \ref{sec2}, we introduce  the {two} models for {SVM}, i.e., L1-loss SVC and $\epsilon$-L1-loss SVR, and give some preliminaries. In Section \ref{sec3}, we apply the {augmented Lagrangian method} to solve {the L1-loss SVC model}. In Section \ref{sec4}, we discuss the semismooth Newton method for subproblem as well as the computational complexity and convergence rate. In Section \ref{sec-svr}, we apply the above framework to the  {$\epsilon$-L1-loss SVR model}.  Numerical results are reported in Section \ref{sec5} to show the efficiency of the proposed method. Final conclusions are given in Section \ref{sec6}.

{\bf Notations.} We use $\|\cdot\|$ to denote the $l_2$ norm for vectors and Frobenius norm for matrices. $\|x\|_1$ is the $l_1$ norm for vector $x$ and $\|x\|_\infty$ is the infinite norm of $x$. $|\Omega|$ denotes the number of elements in set $\Omega$ and $|a|$ denotes the absolute value of the real number $a$. Let $\mathcal S^n$ be the set of symmetric matrices. We use $A\succeq 0$ ($A\succ 0$) to mean that $A\in\mathcal S^n$ is positive semidefinite (positive definite). Let $\hbox{Diag}(u) $ denote a diagonal matrix with diagonal elements coming from vector $u\in\Rm^m$. We use $p^*(\cdot)$ to denote the Fenchel conjugate of a function $p(\cdot)$. 
\section{Problem statement and preliminaries}\label{sec2}
In this section, we will briefly describe two models for SVM and give some preliminaries including Danskin theorem, semismoothness and proximal mapping.

\subsection{Two models for SVM}\label{sec2-1}

\subsubsection{The L1-loss SVC model}

Given training data $(x_i, y_i), \ i = 1, \dots m$, where $x_i
\in\Rm^n$ are the observations, $y_i\in\{-1,1\}$ are the labels, the support vector classification is to find a hyperplane $y = w^Tx + b$ such that the data with different labels can be separated by the hyperplane. The typical SVC model is
\be\label{svm}
\begin{split}
	\min_{\omega\in\Rm^n,\ b\in \Rm} &\ \  \frac12\|\omega\|^2  \\
	\hbox{s.t.}\ \ \ \ \ &  y_i(\omega^Tx_i+b)\ge 1,\ i= 1,\dots,\ m.
\end{split}
\ee
Model (\ref{svm}) is based on the assumption that the two types of data can be successfully separated by the hyperplane. However, in practice, this is usually not the case. A more practical and popular model is the regularized penalty model
\be\label{hing}
\min_{\omega\in\Rm^n,\ b\in\Rm} \\ \  \frac12\|\omega\|^2 + C\sum_{i = 1}^{m} \xi(\omega; x_i, y_i, b),
\ee
where $C>0$ is a penalty parameter and $\xi(\cdot)$ is the loss function. Three frequently used loss functions are as follows.
\bit
\item {L1-loss or $l_1$ hinge loss}:  $\xi(\omega; x_i, y_i, b)=\max(1-y_i(\omega^Tx_i+b), 0)$;

\item {L2-loss or squared hinge loss}:  $\xi(\omega; x_i, y_i, b)=\max(1-y_i(\omega^Tx_i+b), 0)^2 $;
\item {Logistic loss}:  $\xi(\omega; x_i, y_i, b)=\log(1+e^{-y_i(\omega^Tx_i+b)})$.
\eit
As mentioned in Introduction, the L2-loss model has already been solved by semismooth Newton's method in \cite{YinLi2019}.
In our paper, we focus on the L1-loss SVC model, i.e.,
\be\label{prob-svc}
\min_{\omega\in\Rm^n, b\in\Rm} \ \frac12\|\omega\|^2 + C\sum_{i = 1}^m \max(1-y_i(\omega^Tx_i+b),0).
\ee
Notice that there is a bias term $b$ in the standard SVC model. For large-scale SVC, the bias term is often omitted \cite{DCDSVR,LinLiblinear}. By setting
\[
x_i \longleftarrow [x_i,1],\ \  \omega\longleftarrow [\omega,b],
\]
we reach the following model (referred as L1-Loss SVC {\cite{LinLiblinear}})
\be\label{prob}
\min_{\omega\in\Rm^n} \frac12\|\omega\|^2+C\sum_{i=1}^{m} \max(1-y_i(\omega^Tx_i),0).\\
\ee

Mathematically, there is significant difference between problem (\ref{prob-svc}) and problem (\ref{prob}), since problem (\ref{prob-svc}) is convex and problem (\ref{prob}) is strongly convex. On the other hand, it is shown in \cite{DCDSVR} that the bias term hardly affects the performance in most data (See section 4.5 in \cite{DCDSVR} for the numerical comparison with and without bias term). As a result, in our paper, we will focus on the unbiased model (\ref{prob}), which enjoys nice theoretical properties.

\subsubsection{The $\epsilon$-L1-loss SVR model}

Given training data $(x_i,y_i),\ i = 1,\dots, m$, where $x_i\in\Rm^n$, $y_i\in\Rm$, SVR is to find $\omega\in\Rm^n$ and $b\in\Rm$ such that $\omega^Tx_i+b$ is close to the target value $y_i$, $i = 1,\dots, m$. The $\epsilon$-L1-loss SVR model \footnote{For the same reason as in the L1-Loss SVC model, we omit the bias term $b$.} is as follows
\be\label{prob-svr}
\min_{\omega\in\Rm^n} \frac12\|\omega\|^2+C\sum_{i=1}^{m} \max(|\omega^Tx_i-y_i|-\epsilon,0),\\
\ee
where $\epsilon>0$ and $C>0$ are given parameters. We refer to (\ref{prob-svr}) as $\epsilon$-L1-loss SVR as in \cite{yj14}.

\subsection{Preliminaries} \label{sec2-3}
\subsubsection{Semismoothness}

The concept of semismoothness was introduced by Mifflin \cite{Mifflin} for functionals. It was extended to vector-valued functions by Qi and Sun \cite{QiSun93}. Let $\mathcal X$ and $\mathcal Y$ be two real finite dimensional Euclidean spaces with an inner product $\langle\cdot,\cdot\rangle$ and its induced norm $\|\cdot\|$ on $\mathcal X$.

\begin{definition}
	(Semismoothness \cite[]{Mifflin,QiSun93,semimatSS}). Let $\Phi:\mathcal O\subset\mathcal X\to\mathcal Y$ be a locally Lipschitz continuous function on the open set $\mathcal O$. We say that $\Phi$ is semismooth at $x$ if (i) $\Phi$ is directional differentiable at $x$ and (ii) for any $V\in\partial \Phi(x+h)$,
	\[
	\Phi(x+h) -\Phi(x) -Vh = o(\|h\|),\ h\to0.
	\]
Here $\partial \Phi(x+h)$ is the Clarke subdifferential  \cite{Clarke} of $\Phi$ at $x+h$.
	$\Phi$ is said to be strongly semismooth at $x$ if $\Phi$ is semismooth at $x$ and for any $V\in\partial \Phi(x+h)$,
	\[
	\Phi(x+h) -\Phi(x) -Vh = O(\|h\|^2),\ h\to0.
	\]
\end{definition}

It is easy to check that piecewise linear functions are strongly semismooth. Furthermore, the composition of (strongly) semismooth functions is also (strongly) semismooth.
A typical example of strongly semismooth function is  $\max (0, t),\ t\in\Rm$.

\subsubsection{Moreau-Yosida regularization}

\begin{figure}[htbp]
	\centering
	\begin{minipage}[t]{0.48\textwidth}
		\centering
		\includegraphics[width=8cm]{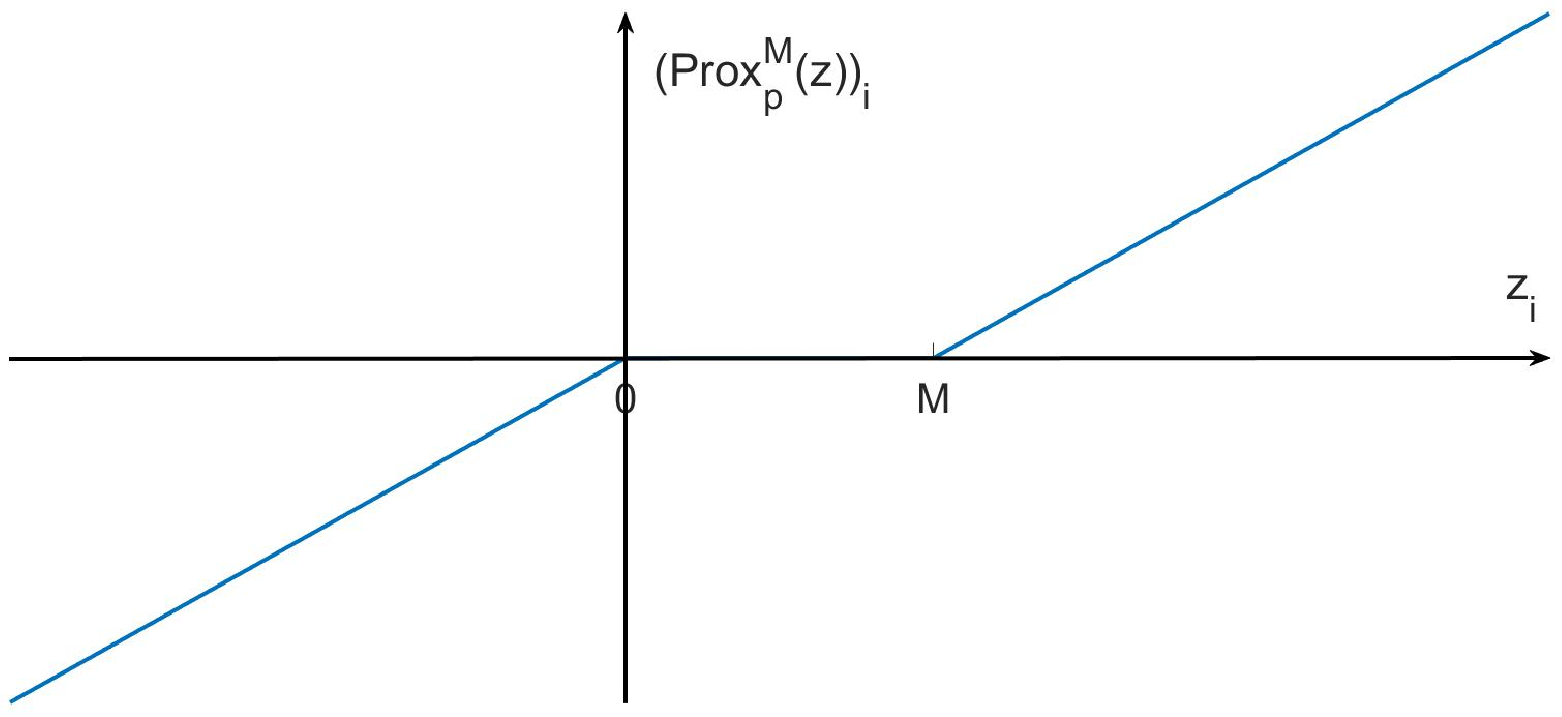}
		\caption{Demonstration of $\hbox{Prox}_p^M(z_i)$ }
		\label{Prox-svc-fig}
	\end{minipage}
	\begin{minipage}[t]{0.48\textwidth}
		\centering
		\includegraphics[width=8cm]{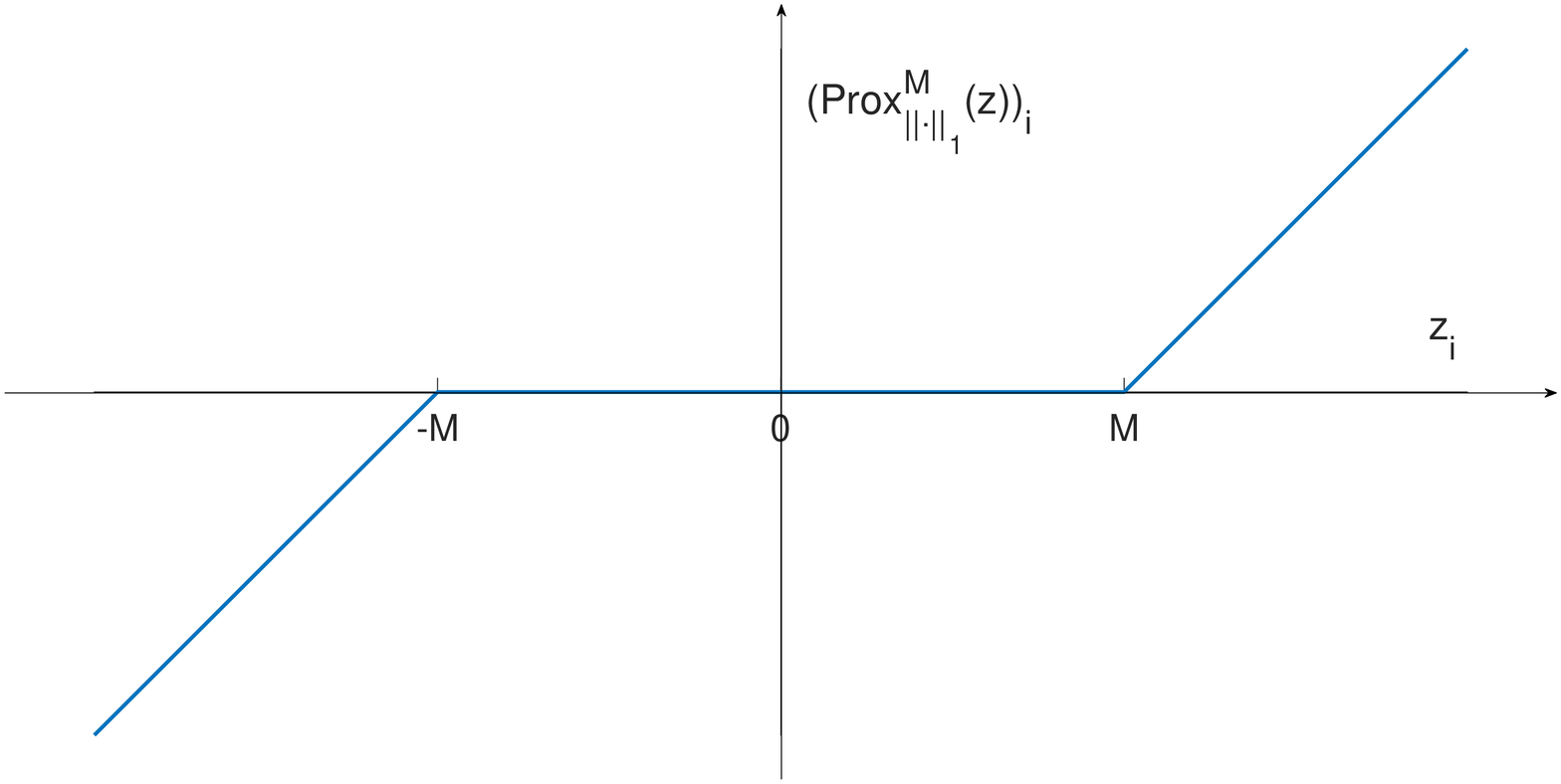}
		\caption{Demonstration of $\hbox{Prox}_{\|\cdot\|_1}^M(z_i)$}
		\label{prox-l1norm}
	\end{minipage}
\end{figure}

Let $q:\mathcal X\longrightarrow (-\infty,+\infty)$ be a closed convex function. The Moreau-Yosida \cite{moreau,yosida} regularization of $q$ at $x\in\mathcal X$  is defined by
\be\label{MY-theta}
\theta_q(x):=\min_{y\in\mathcal X} \frac12\|y-x\|^2+q(y).
\ee
The unique solution of (\ref{MY-theta}), denoted as $\hbox{Prox}_q(x)$, is called the proximal point of $x$ associated with $q$. The following property holds for Moreau-Yosida regularization \cite[Proposition 2.1]{ppaNNM2009}.
\begin{proposition}\label{prop-MY}
Let $q:\mathcal X\longrightarrow (-\infty,+\infty)$ be a closed convex function, $\theta(\cdot)$ be the Moreau-Yosida regularization of $q$ and $\hbox{Prox}_q(\cdot)$  be the associated proximal point mapping. Then $\theta_q(\cdot)$ is continuously differentiable, and there is
\[
\nabla\theta_q(x) = x-\hbox{Prox}_q(x).
\]
\end{proposition}
Let
	\be\label{p}
	p(s) = C\sum_{i = 1}^m\max(s_i,0), \ s\in \Rm^m.
	\ee
The proximal mapping, denoted as $\hbox{Prox}_p^M(\cdot)$,  is defined as  the solution of the following problem

\be\label{prob-Prox}
\phi(z):=\min_{s\in\Rm^m}\psi(z,s), 
\ee
where $\psi(z,s):=\frac1{2M}\|z-s\|^2 +  p(s),\ z\in\Rm^m.$ It is easy to derive that $\hbox{Prox}_p^M(z)$ takes the following form ({See Appendix for the details of  deriving $\hbox{Prox}_p^M(z)$})
\be\label{Prox-svc}
(\hbox{Prox}_p^M(z))_i=\left\{
\begin{array}{ll}
	z_i-CM,\ & z_i>CM,\\
	z_i,\ & z_i<0,\\
	0,\ & 0\leq z_i\leq CM.
\end{array}
\right.
\ee
{The proximal mapping $\hbox{Prox}_p^M(\cdot)$} is  piecewise linear as shown in  Figure \ref{Prox-svc-fig}, and therefore strongly semismooth. 
\begin{remark}Here we would like to point out that the  proximal mapping here is closely related to  the proximal mapping associate with  $\|s\|_1$, which is popular used in the well-known LASSO problem. Actually, from Figure \ref{Prox-svc-fig} and Figure \ref{prox-l1norm}, the proximal mapping $\hbox{Prox}_p^M(\cdot)$ can be viewed as a shift of that with $\|s\|_1$.
\end{remark}

Given $z$, the Clarke subdifferential of $\hbox{Prox}_p^M(z)$, denoted as $\partial \hbox{Prox}_p^M(z)$, is a set of diagonal matrices. For $U\in \partial \hbox{Prox}_p^M(z)$, its diagonal elements take the following form
	\be
	U_{ii}=\left\{
	\begin{array}{ll}
		1,\ & z_i>CM, \hbox{ or }z_i<0,\\
		0,\ & 0< z_i<CM, \\
		{ u_i,\ 0\le u_i\le1},\ & z_i=0, \hbox{ or } z_i=CM.
	\end{array}
	\right.
	\ee

In other words, we have
\begin{eqnarray*}
	\partial \hbox{Prox}_p^M(z)& =& \big\{\hbox{Diag}(u): \ u\in \Rm^m, u_i = 1, \hbox{ if } z_i>CM \hbox{ or } z_i<0; \\
	&& {u_i = 0, \hbox{ if } 0<z_i<CM; \ u_i\in[0,1], \ \hbox{otherwise}} \big\}.
\end{eqnarray*}

\begin{figure}[htbp]
	\centering
	\includegraphics[width=8cm]{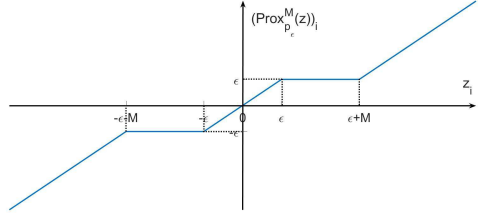}
	\caption{Demonstration of ${\rm Prox} _{p_{\epsilon}}^M(z_i)$}
	\label{Prox-svr-fig}
\end{figure}

Similarly, for $p_{\epsilon}(\cdot)$ defined by
\be\label{p-e}
p_{\epsilon}(s) = C\sum_{i =1}^m\max(0,|s_i|-\epsilon),
\ee
there is (see Figure \ref{Prox-svr-fig} for
${\rm Prox} _{p_{\epsilon}}^M(z)$)
\[
({\rm Prox} _{p_{\epsilon}}^M(z))_i=\left\{
\begin{array}{ll}
z_i-CM,&z_i\ge\epsilon+CM,\\
\epsilon,&\epsilon< z_i<\epsilon+CM ,\\
z_i,&| z_i|\le\epsilon,\\
-\epsilon,&-\epsilon-CM<z_i< -\epsilon,\\
z_i+M,& z_i\le-\epsilon-CM.
\end{array}
\right.
\]It can be easily verified that $\partial {\rm Prox} ^M_{p_{\epsilon}}(z)$ is strongly semismooth as well. For $U\in\partial {\rm Prox} ^M_{p_{\epsilon}}(z)$, its diagonal elements are given by
\[
U_{ii}=\left\{
\begin{array}{ll}
1,\ & |z_i|>\epsilon+CM, {\hbox{ or }|z_i|<\epsilon},\\
0,\ & \epsilon< |z_i|<\epsilon+CM,\\
{ u_i,\ 0\le u_i\le 1,}\ & |z_i|=\epsilon+CM,\ \hbox{ or } |z_i|=\epsilon.\\
\end{array}
\right.
\]



We end this section by the the following property of $p(\cdot)$.
\begin{proposition}\label{lem-1}
Let \[P=\{x\in\Rm^m\ |\ \|x-\frac C2 e\|_\infty\le \frac C2\},\] where $e=(1,\dots, 1)^T\in\Rm^m$. Then we have \[p(\cdot)=\delta_P^*(\cdot).\]
\end{proposition}

{\bf Proof.} Note that $p(\cdot)$ can be equivalently written as
\[
p(s) = \frac C2 s^Te + \frac C2 \|s\|_1.
\]
By the definition of Fenchel conjugate function, there is

\begin{eqnarray*}
p^*(y)&=& \sup_{s\in\Rm^m}  s^Ty - p(s) \\
&=& \sup_{s\in\Rm^m}  s^Ty- \frac C2 s^Te - \frac C2 \|s\|_1 \\
&=&\sup_{s\in\Rm^m}  s^T(y- \frac C2 e) - \frac C2 \|s\|_1 \\
&=&\frac C2 \sup_{s\in\Rm^m}  \frac 2C s^T(y- \frac C2 e) -  \|s\|_1 \\
&:=& \frac C2 \sup_{s\in\Rm^m}  \frac 2C s^T(y- \frac C2 e) -  u(s)  \quad\ (u(s):=\|s\|_1)\\
&=&\frac C2 u^*(\frac 2C(y- \frac C2 e)).
\end{eqnarray*}

Note that $u^*(\cdot) = \delta_{P_0}(\cdot)$, where $P_0:=\{x\in\Rm^m\ |\ \|x\|_\infty \le 1\}$, and $\delta_{P_0}(\cdot)$  is the indicator function defined as $0$ if $x\in P_0$ and $+\infty$ otherwise. We then have
\[p^*(y) =\frac C2 \delta_{P_0}(\frac 2C(y- \frac C2 e))= \delta_{P_0}(\frac 2C(y- \frac C2 e))=\delta_P(y).\]
In other words, there is
$p(\cdot)=\delta_P^*(\cdot)$. The proof is finished. \hfill$\Box$

\begin{proposition}\label{p_conjugate_svr}
Let $\overline{P}:=\{x\in\Rm^m:\|x\|_{\infty}\le C \}$. There is
\[p_{\epsilon}^*(y)
=\left\{\begin{array}{ll}{\epsilon\|y\|_1,} & {y\in \overline{P},} \\
{+\infty,} & {\text{otherwise}.}\end{array}\right.
\]
\end{proposition}

{\bf Proof.} Similar to the proof of Proposition \ref{lem-1}, there is
\begin{eqnarray*}
	p_{\epsilon}^*(y)&=& \sup _{x \in \Rm^{m}} x^Ty-C\sum_{i=1}^{m}\max(0,|x_i|-\epsilon) \\
	&=& \sum_{i=1}^{m} \sup _{x_{i} \in \Rm} x_{i} y_{i}-C \max (0,|x_{i}|-\epsilon).
\end{eqnarray*}
Note that
\[x_iy_i-C\max(0,|x_i|-\epsilon)
=\left\{\begin{array}{ll}{x_{i} y_{i},} & {\left|x_{i}\right| \leq \epsilon,} \\ {x_{i}\left(y_{i}-C\right)+C \epsilon,} & {x_{i}>\epsilon,} \\ {x_{i}\left(y_{i}+C\right)+C \epsilon,} & {x_{i}<-\epsilon.}\end{array}\right.
\]
We consider the following different cases of $y_i$:

\bit
\item if $|y_i|=C$, there is $\sup_{x_i} x_{i} y_{i}-C\max(0,|x_{i}|-\epsilon)=C\epsilon$;

\item if $|y_i|<C$, there is $\sup_{x_i} x_{i} y_{i}-C\max(0,|x_{i}|-\epsilon)=|y_i|\epsilon$;

\item if $|y_i|>C$, $\sup_{x_i} x_{i} y_{i}-C\max(0,|x_{i}|-\epsilon)=+\infty$.
\eit

To sum up, we have
\[\sup_{x_i\in\Rm} x_iy_i-C\max(0,|x_i|-\epsilon)
=\left\{\begin{array}{ll}{|y_i|\epsilon,} & {|y_i|\le C,} \\
{+\infty,} & {|y_i|>C.} \end{array}\right.
\]
Consequently, we get
\[p_{\epsilon}^*(y)
=\left\{\begin{array}{ll}{\epsilon\|y\|_1,} & {\|y\|_{\infty}\le C,} \\
{+\infty,} & {\text{otherwise}.} \end{array}\right.
\]

The proof is finished. \hfill$\Box$


\section{The augmented Lagrangian method}\label{sec3}
In this section, we will discuss the augmented Lagrangian method to solve the L1-loss SVC model (\ref{prob}).

\subsection{Problem reformulation}
{We first reformulate  (\ref{prob}) equivalently as the following form}
\be\label{prob-1}
\min_{w\in \Rm^n}\frac{1}{2}\|w\|^2+C\|\max (0,Bw+d)\|_1
\ee
by letting
$$
B=-\left[
\begin{array}{c} 
y_1x_1^T  \\
\vdots \\
y_mx_m^T
\end{array}
\right]\in\Rm^{m\times n}
,\
d=\left[
\begin{array}{c}
1  \\
\vdots \\
1
\end{array}
\right]\in\Rm^m.
$$
By introducing a variable $s\in\Rm^m$, we get the following {constrained} optimization problem
\be\label{prob-2}
\begin{aligned}\min_{w\in \Rm^n,{s\in\Rm^m}}&\  \frac{1}{2}\|w\|^2+p(s)
	\\
	\hbox{s.t.}&\ s=Bw+d.
\end{aligned}
\ee
where $p(s)$ is defined as in (\ref{p}). The Lagrangian function for (\ref{prob-2}) is
\[
l(w,s,\lambda) = \frac12\|w\|^2+p(s) -\langle \lambda, s-Bw-d\rangle,
\]
where $\lambda\in\Rm^m$ is the Lagrange multiplier corresponding to the equality constraints.
The dual problem is
\be\label{dual}
\max_{\lambda\in\Rm^m}\ -f(\lambda):=-\frac{1}{2}\|B^T\lambda\|^2+\left\langle \lambda,d \right\rangle -p^*(\lambda).
\ee
The KKT conditions associated with problem (\ref{prob-2}) are given by
\be\label{KKT}
\left\{
\begin{array}{ll}
s=Bw+d,\\
w+B^T\lambda=0,\\
0\in\partial p(s)-\lambda.
\end{array}
\right.
\ee

\begin{remark}\label{remark3-1}
It is easy to see that both (\ref{prob-2}) and (\ref{dual}) admit feasible solutions. Consequently, both (\ref{prob-2}) and (\ref{dual}) admit optimal solutions. By Theorem 2.1.8 in \cite{Sun1986}, there is no duality gap between (\ref{prob-2}) and (\ref{dual}). Moreover, by Theorem 2.1.7 in \cite{Sun1986}, $(s^*,w^*,\lambda^*)$ solves (\ref{KKT}) if and only if $(s^*,w^*)$ is an optimal solution to (\ref{prob-2}) and $\lambda^*$ is an optimal solution of (\ref{dual}). Consequently, the set of Lagrange multipliers is not empty.
\end{remark}

\begin{remark}\label{remark-1}
We would like to point out that after we finished our paper, we realized that problem (\ref{prob-2}) is a special case of the general problem considered in \cite{LiSun2017}. Consequently, as we will show below, it enjoys interesting properties, due to  which  the convergence results of the augmented Lagrangian method (Theorems \ref{thm-global} and \ref{thm-local}) in our paper directly follows from those in \cite{LiSun2017}.
\end{remark}

\begin{proposition}\label{prop-soc}
Let $(w^*, s^*)$ be the optimal solution to problem (\ref{prob-2}).  Then the second order sufficient condition holds at {$(w^*, s^*)$}.
\end{proposition}
{\bf Proof.} Assume that $(w^*, s^*, \lambda^*)$ is one of the solutions for the KKT system (\ref{KKT}). Note that the effective domain of $p(\cdot)$, denoted as $\hbox{dom}(p)$,  is $\Rm^m$, therefore,  the tangent cone of $\hbox{dom}(p)$ at $s^*$ is $\Rm^m$. That is $T_{\hbox{dom}(p)}(s^*)=\Rm^m$. As a result, by the definition in \cite[(7)]{LiSun2017},
the critical cone associate with (\ref{prob-2}) at $(w^*,s^*)$, denoted by $\mathcal{C}(w^*,s^*)$, is as follows
\[
\mathcal{C}(w^*,s^*) =\left\{
(d_1,d_2)\in\Rm^n\times \Rm^m\ | \ Bd_1-d_2=0, \ (w^*)^Td_1+p'(s^*;d_2)=0
\right\}.
\]
Here, $p'(s^*;d_2)$ denote the directional derivative of $p$ at $s^*$ with respect to $d_2$.

Next, we will show that
\be\label{eq-4}
(d_1,d_2)\neq0\Longrightarrow d_1\neq0, \ \forall \ (d_1,d_2)\in \mathcal{C}(w^*,s^*).
\ee
For contradiction, assume that for any $(d_1,d_2)\in \mathcal{C}(w^*,s^*)$, $d_1=0$. Then $(d_1,d_2)\in \mathcal{C}(w^*,s^*)$ gives that $d_2 = Bd_1=0$, contradicting with $(d_1,d_2)\neq0$. Consequently, (\ref{eq-4}) holds.

By the definition in \cite[Definition 2.6]{LiSun2017}, the second order sufficient condition holds at $(w^*,s^*)$ if
\be\label{ineq-1}
\langle d_1, (\nabla h(w^*))' d_1\rangle >0, \ \forall\ 0\neq (d_1,d_2)\in\mathcal{C}(w^*,s^*), \ \hbox{where} \ h(w)=\frac12\|w\|^2.
\ee
By the definition of directional derivative, it is easy to calculate that \[
(\nabla h(w^*))' d_1=\lim_{t\downarrow 0}\frac{w^*+td_1- w^*}{t}=d_1.
\]
Together with (\ref{eq-4}), (\ref{ineq-1}) reduces to
\[
\langle d_1,d_1\rangle>0, \forall \ d_1\neq 0, (d_1,d_2)\in \mathcal{C}(w^*,s^*),
\]
which  obviously holds. The proof is finished. \hfill$\Box$

\subsection{Augmented Lagrangian method (ALM)}
Next, we will apply {augmented Lagrangian method (ALM)} to solve (\ref{prob-2}).
The augmented Lagrangian function of (\ref{prob-2}) is
{
\[
L_{\sigma}(w,s;\lambda)=\frac{1}{2}\|w\|^2+p(s)-\langle\lambda,s-Bw-d\rangle+\frac{\sigma}{2}\|s-Bw-d\|^2,
\]
}
where {$\sigma>0$}.

ALM  works as follows. At iteration $k$,  solve
\be\label{prob-sub}
\min_{w,s} \ L_{\sigma_k}(w,s;\lambda^k)
\ee
to get $(w^{k+1},s^{k+1})$. Then update the Lagrange multiplier by \[\lambda^{k+1}=\lambda^k-\sigma_k(s^{k+1}-Bw^{k+1}-d),\] { and $\sigma_{k+1}\ge\sigma_{k}$}.

The key step in ALM is to solve the subproblem (\ref{prob-sub}). Similar to that in \cite{LiSun2017},  given fixed $\sigma>0$ and $\lambda$, let $(w^*,s^*)$ denote the unique solution of subproblem (\ref{prob-sub}). Denote
 \[
 \phi(w):=\min_sL_{\sigma}(w,s;\lambda).
 \] Note that \[\min_s\  L_{\sigma}(w,s;\lambda)=\frac12\|w\|^2-\frac{1}{2\sigma}\|\lambda\|^2+\sigma\min_s\ \frac{1}{2}\|s-z(w)\|^2+\frac{1}{\sigma}p(s),\]
where \[
z(w)=Bw+d+\frac{\lambda}{\sigma}.
\]
There is
\be\label{phi}
\phi(w)=\frac{1}{2}\|w\|^2-\frac{1}{2\sigma}\|\lambda\|^2+\sigma\tau(z(w)),
\ee
where $\tau(z(w)):=\min_s\ \frac{1}{2}\|s-z(w)\|^2+\frac{1}{\sigma}p(s)$ is the Moreau-Yosida regularization of $\frac{1}{\sigma}p(s)$. Therefore, we can get $(w^*,s^*)$ by
\begin{align}
w^*&=\arg\min\phi(w),\label{sub}\\
s^*&=s^*(w^*)=\hbox{Prox}^{1/\sigma}_p(z(w^*)).\label{s}
\end{align}

Due to the strong convexity and the differentiability of $\phi(w)$ in Proposition \ref{prop-MY}, we will apply semismooth Newton's method  to solve (\ref{sub})
in Section \ref{sec4}.  {The details  of ALM is given in Algorithm \ref{alg-alm-1}.}

\begin{algorithm}\label{alg-alm-1} { ALM for L1-loss SVC}
	\bit
	
	\item [S0] Initialization. 
$\sigma_0>0,\lambda^0\in\Rm^m;w^0\in\Rm^n.$ {For $k_1=0, 1, 2, ...$}
	
	{
	\item [S1] Solve (\ref{sub}) to get $w^{k+1}$. Then Calculate $s^{k+1}=s^*(w^{k+1})$ by (\ref{s}).
	}
	%
	\item [S2] Update $\lambda^{k+1}=\lambda^k-\sigma_k(s^{k+1}-Bw^{k+1}-d)$.
	
	\item [S3] {Choose $\sigma_{k+1}$ such that $\sigma_{k}\le\sigma_{k+1}<+\infty.$  Go to S1.}
%
	\eit 
\end{algorithm}

\subsection{Convergence results}
To guarantee the global convergence of ALM, the following standard stopping criteria \cite{Roch1976a,Roch1976b} are used in \cite{LiSun2017} to solve (\ref{prob-2}) approximately
\[
(A)\ \ \Psi_k(w^{k+1},s^{k+1}) -\inf\Psi_k \le \epsilon^2_k/ 2\sigma_k, \ \sum_{k=0}^\infty\epsilon_k<+\infty
\]
where $\Psi_k(w,s):= L_{\sigma_k}(w,s,\lambda^k)$.
The global convergence result of ALM was originally from \cite{Roch1976a, Roch1976b}. As we mentioned in Remark \ref{remark-1}, ALM is   applied to solve a general form of (\ref{prob-2}), i.e., problem (D) in \cite{LiSun2017}, and the convergence result, i.e., Theorem 3.2 in \cite{LiSun2017} is obtained therein. Therefore, as a special case of the problem (D) in \cite{LiSun2017},  we can get the  following global convergence result of Alg. \ref{alg-alm-1}.

\begin{theorem}\label{thm-global}
Let the sequence $\{(w^k,s^k, \lambda^k)\}$ be  generated by ALM with stopping criteria (A). Then the sequence $\{\lambda^k\}$  is bounded and converges to an optimal solution of (\ref{dual}). Moreover, the sequence $\{(w^k,s^k)\}$ is also bounded and converges to the unique optimal solution $(w^*,s^*)$ of (\ref{prob-2}).
\end{theorem}
{\bf Proof.} Note that the optimal solution for strongly convex problem (\ref{prob-1}) exists and it is unique. Consequently, (\ref{prob-2}) admits a unique solution. By Remark \ref{remark3-1}, the set of Lagrange multiplier is also non-empty. Consequently, by Theorem 3.2 in \cite{LiSun2017}, the result holds. \hfill $\Box$

To state the local convergence rate, we need the following stopping criteria which are popular used such as in \cite{Sun2018} and \cite{LiSun2017}.
\[
(B1)\ \ \Psi_k(w^{k+1},s^{k+1})-\inf\Psi_k\le (\delta_k^2/2\sigma_k)\|\lambda^{k+1}-\lambda^k\|^2, \ \sum_{k=1}^\infty\delta_k<+\infty,
\]
\[
(B2) \ \ \hbox{dist}(0,\partial\Psi_k(w^{k+1},s^{k+1}))\le(\delta'_k/\sigma_k)\|\lambda^{k+1}-\lambda^k\|, \ \sum_{k=1}^\infty\delta_k'\rightarrow 0.
\]
We also need the following definitions.
Define the maximal monotone operator $T_f$ and $T_l$ \cite{Roch1976a} by
\[
T_f = \partial f(\lambda), \ T_l (w,s,\lambda) = \{(w',s',\lambda')\ | \ (w',s',-\lambda')\in\partial l(w,s,\lambda)\}.
\]
\begin{definition}
Let $F:\mathcal X\rightrightarrows \mathcal Y$ be a multivalued mapping and $(\widetilde x,\widetilde y)\in \hbox{gph} F$, where $\hbox{gph} F$ is the graph of $F$ defined by
\[
\hbox{gph} F:=\{(x,y)\in\mathcal X\times \mathcal Y\ |\ y\in F(x)\}.
\]
$F$ is said to be metrically subregular at $\widetilde x$ for $\widetilde y$ with modulus $\kappa\ge0$ if there exist neighbourhoods $U$ of $\widetilde x$ and $V$ of $\widetilde y$ such that
\[
\hbox{dist}(x,F^{-1}(\widetilde y)) \le\kappa \hbox{dist}(\widetilde y, F(x)\cap V), \ \forall \ x\in U.
\]
\end{definition}
\begin{definition}
Let $F:\mathcal X\rightrightarrows\mathcal Y$ be a multivalued mapping and $y\in\mathcal Y$ satisfy $F^{-1}(y)\neq \emptyset$. $F$ is said to satisfy the error bound condition for the point $y$ with modulus $\kappa\ge 0$ if there exists $\epsilon>0$ such that  if $x\in\mathcal X$ with $\hbox{dist}(y, F(y))\le\epsilon$, then
\[
\hbox{dist}(x,F^{-1}(y))\le \kappa \hbox{dist}(y,F(x)).
\]

\end{definition}

\begin{proposition}\label{prop-1}
Let $F:\mathcal X\rightrightarrows\mathcal Y$ be a ployhedral multifunction. Then $F$ satisfies the error bound condition for any $y\in \mathcal Y$ satisfying $F^{-1}(y)\neq \emptyset$ with a common modulus $\kappa\ge0$.
\end{proposition}

Now we are ready to give the convergence rate of ALM, which is similar to that in \cite[Theorem 3.3]{LiSun2017}.
\begin{theorem}\label{thm-local}
Let $\{(w^k,s^k, \lambda^k)\}$ be the infinite sequence generated by ALM with stopping criteria (A) and (B1). The following results hold.
 \bit
 \item [(i)] The sequence $\{\lambda^k\}$   converges to $\lambda^*$, one of the  optimal solutions of (\ref{dual}), and for all $k$ sufficiently large, there is
 \[
 \hbox{dist}(\lambda^{k+1},\Omega)\le\beta_k\hbox{dist}(\lambda^k,\Omega),
 \]
where $\Omega$ is the set of the optimal solutions of (\ref{dual}) and \[\beta_k=\frac{a_f(a_f^2+\sigma_k^2)^{-\frac12}+2\delta_k}{1-\delta_k}\to \beta_\infty=a_f(a_f^2+\sigma_\infty^2)^{-\frac12}<1\]
as $k\to+\infty$.
 Moreover, $\{(w^k,s^k)\}$  converges to the unique optimal solution $(w^*,s^*)$ of (\ref{prob-2}).
\item [(ii)]
If the stopping criteria (B2) is also used, then for all $k$ sufficiently large, there is
\[
\|(w^{k+1},z^{k+1})-(w^*,z^*)\|\le\beta'_k\|\lambda^{k+1}-\lambda^k\|,
\]
where $\beta'_k=a_l(1+\delta'_k)/\sigma_k \to  a_l/\sigma_\infty$ as $k\to +\infty$.
\eit
\end{theorem}

{\bf Proof.} (i) We only need to show that $T_f$ satisfies the error bound condition for the origin with modulus $a_f$. By Definition 2.2.1 in \cite{Sun1986}, (\ref{dual}) is  convex piecewise linear-quadratic programming problem.  By \cite[Proposition 2.2.4]{Sun1986}, we know that the corresponding operator  $T_f$ is polyhedral multivalued function. Therefore, by Proposition \ref{prop-1}, the error bound condition holds at the origin with modulus  $a_f$. Following Theorem 3.3 in \cite{LiSun2017}, we proved (i).  To show (ii), first note that the second order sufficient condition holds for problem (\ref{prob-2}) by Proposition \ref{prop-soc}. On the other hand, by Proposition \ref{lem-1}, $p(\cdot)$ is the support function of a non-empty polyhedral convex set. Consequently, by Theorem 2.7 in \cite{LiSun2017}, $T_l$ is metrically subregular at $(w^*,s^*, \lambda^*)$ for the origin. Then by the second part of  Theorem 3.3 in \cite{LiSun2017}, (ii) holds. The proof is finished. \hfill $\Box$


\section{Semismooth Newton's method for solving  (\ref{sub})}\label{sec4}
In this section, we will discuss semismooth {Newton's} method to solve (\ref{sub}). In the first part, we will give the details of semismooth {Newton's method}. In the second part, we will analyse the computational complexity of the key step in semismooth Newton's method. The third part is devoted to the convergence result.
\subsection{ Semismooth Newton's method}
Note that $\phi(w)$ defined as in (\ref{phi}) is continuously differentiable. By Proposition \ref{prop-MY}, the gradient $\nabla \phi(w)$ takes the following form
\be
\nabla\phi(w) = w+\sigma\nabla\tau(z(w))=
w+B^T\lambda-\sigma B^T(s^*(w)-Bw-d),
\ee
which is strongly semismooth due to the strongly semismoothness of $s^*(w)$.
The generalized Hessian of $\phi$ at $\omega$, denoted as $\partial^2\phi(\omega)$, satisfies the following condition \cite[Proposition 2.3.3, Proposition 2.6.6]{Clarke}: for any $h\in \Rm^n$,
\[
\partial^2\phi(\omega)(h)\subseteq \hat\partial^2\phi(\omega)(h),
\]
where
\be\label{v}
\hat\partial^2\phi(\omega) = I + \sigma B^TB-\sigma B^T\partial \hbox{Prox}^{1/\sigma}_p(z(w))B.
\ee
Consequently, solving (\ref{sub}) is equivalent to solving the following strongly semismooth equations
\be\label{sub-w}
\nabla \phi(w)=0.
\ee

We apply semismooth Newton's method to solve the non-smooth equations (\ref{sub-w}). At iteration $k$,  we update $w$ by
\[
w^{k+1} = w^k -(V^k)^{-1} \nabla \phi(w^k), \ V^k\in\hat\partial^2 \phi(w^k).
\]
Below,
we use the following
well studied globalized version of the semismooth Newton method \cite[Algorithm 5.1]{QiSun06} to solve (\ref{sub-w}). 

\begin{algorithm}\label{algorithm2}  A globalized semismooth Newton method
	\begin{itemize}
		
		\item [S0] Given $j:=0$. Choose $\omega^0,\ \sigma\in(0,\ 1)$, $\rho\in(0,\ 1)$,\ $\delta>0$, $\eta_0>0$ and $\ \eta_1>0 $.
		\item [S1] Calculate $\nabla \phi(\omega^j)$. If $\|\nabla \phi(\omega^j)\|\le \delta$, stop. Otherwise, go to S2.
		\item [S2] Select an element $V^j\in \hat\partial^2\phi(\omega^j)$ defined as in (\ref{v}). Apply conjugate gradient (CG) method \cite{Hestenes} to find an approximate solution $d^j$ by
		\be\label{linearsystem}
		V^j d^j + \nabla \phi(\omega^j)=0
		\ee
		such that
		\[
		\|V^j d^j + \nabla \phi(\omega^j)\|\le \mu_j \|\nabla \phi(\omega^j)\|,
		\]
		where  $\mu_j = \min(\eta_0    ,\ \eta_1\|\nabla \phi(\omega^j)\|)$.
		
		\item [S3] Do line search, and let $m_j>0$ be the smallest integer such that the following holds
		\[
		\phi(\omega^j+\rho^md^j)\le \phi(\omega^j) + \sigma \rho^m\nabla \phi(\omega^j)^Td^j.
		\]
		Let $\alpha_j = \rho^{m_j}$.
		\item [S4] Let $\omega^{j+1} = \omega^j + \alpha_j d^j$, $j:=j+1$, go to S1.
		
	\end{itemize}
\end{algorithm}

\subsection{Computational complexity}\label{sec4-2}
{It is well-known that} the high computational complexity for semismooth Newton's method comes from solving the linear system in S2 of Alg. \ref{algorithm2}. As designed in Alg. \ref{algorithm2}, we apply the popular solver conjugate gradient method (CG) to solve the linear system iteratively. Recall that $w\in\Rm^n$, $B\in\Rm^{m\times n}$. The {heavy} burden in applying CG is to calculate $Vh$, where  $h\in\Rm^n$ is a given vector. Below, we will mainly analyse several ways to implement $Vh$, where
\[
V = I + \sigma B^TB - \sigma B^TUB,\quad U\in\partial \hbox{Prox}^{1/\sigma}_p(z(w)).
\]

{\bf Way I.} Calculate $V$ explicitly and then calculate $Vh$.

{\bf Way II.} Do not save $V$ explicitly.  Instead, calculate $Vh$ directly by
\[
Vh = h + \sigma B^T(Bh)- \sigma B^TU(Bh).
\]
Note that $U=\hbox{Diag}(u) $ is a diagonal matrix with vector $u\in\Rm^m$. We can explore the diagonal structure, and use the following formula
\be\label{way2}
Vh = h + \sigma B^T(Bh)- \sigma B^T\hbox{Diag}(u)(Bh).
\ee
Specifically, we can first compute $Bh$, then do the rest computation as in (\ref{way2}).

{\bf Way III.} Do not save $V$. Inspired by the technique in \cite[Section 3.3]{LiSun2017},  we reformulate $Vh$ as
\[
Vh = h+\sigma B^T(I-U)Bh=h + \sigma B^T(I-\hbox{Diag}(u))(Bh).
\]
Reformulate $I-U$ as
\[
\begin{array}{ll}
	I- \hbox{Diag}(u) &= I-\left\{
	\begin{array}{ll}
	1,\ & z_i>\frac{C}\sigma, \hbox{ or }z_i<0,\\
	0,\ & 0<z_i<\frac C\sigma,\\
	u_i,\ u_i\in[0,1],\ & z_i=0, \hbox{ or } z_i = \frac {C}{\sigma},
	\end{array}
	\right.\\
	&=\left\{
	\begin{array}{ll}
	0,\ & z_i>\frac{C}\sigma, \hbox{ or }z_i<0,\\
	1,\ & 0<z_i<\frac C\sigma,\\
	1-u_i\in[0,1],\ & z_i=0, \hbox{ or } z_i = \frac C\sigma.
	\end{array}
	\right.
\end{array}
\]
{By choosing $u_i = 1$ for $z_i=0$ or $z_i=\frac C\sigma$, we get}
\[
I-U = \left\{
\begin{array}{ll}
0,&z_i\ge \frac{C}\sigma, \hbox{ or }z_i\le 0,\\
1,& 0<z_i<\frac C\sigma.
\end{array}
\right.
\]
Now let
\[
\mathcal{I}(z) = \{i: \ z_i\in(0,\frac C\sigma)\}.
\]
Then \[B^T(I-U)Bh=B(\mathcal{I}(z),:)^TB(\mathcal{I}(z),:)h.\]
Finally, we calculate $Vh$ by
\be
Vh = h + \sigma B(\mathcal{I}(z),:)^T(B(\mathcal{I}(z),:)h).
\ee
This lead to the computational complexity {$O(n|\mathcal{I}(z)|)$}.

We summarize the details of the computation for the three ways in Table \ref{tab-cost}. As we will show in numerical part, in most situations, $|\mathcal{I}(z)|$ is far more less than $m$. Together with the complexity reported in Table \ref{tab-cost}, we have the following relations
\[
{O(n|\mathcal{I}(z)|)}<O(nm)< O(nm+m^2)<O(n^2m^2+nm^3).
\]
As a result, we use Way III in our algorithm.
\begin{remark}
We would like to point out that the technique used in Way III essentially follows from the idea in \cite[Section 3.3]{LiSun2017}. The sparse structure is fully explored in semismooth Newton's method to solve the subproblem in ALM \cite{LiSun2017} and it leads to a highly efficient solver for the LASSO  problem.
\end{remark}
\begin{table}
	\centering
	\tbl{Computational Cost for Traditional Implementation}
	{\begin{tabular}{|c|c|c|c|}
		
		\hline
		Way&  Formula & Computational Cost& Complexity
		\\\hline

		Way I&Form {$V=I + \sigma B^TB - \sigma B^TUB$} & $2n^2m^2+n^2m^3$& $O(n^2m^2+nm^3)$\\
		
		&Calculate $V\Delta w${, where $\Delta w=w^{k+1}-w^{k}$}& $n^2$&\\\hline
		
		Way II&Calculate $tmp = Bh$ 
		& $mn$& $O(nm+m^2)$\\
		
		&Calculate $h+\sigma B^T(tmp)-\sigma B^T\hbox{Diag}(v)(tmp)$& $2nm+m^2+n$&\\\hline
		Way III & Calculate $tmp={B(\mathcal{I}(z),:)}h$ & $|\mathcal{I}(z)|n$& $O({|\mathcal{I}(z)|n})$\\
		& Calculate $h+\sigma (B(\mathcal{I}(z),:)^T*tmp)$ &{$|\mathcal{I}(z)|n+n$}&\\
		\hline
	\end{tabular}}
	\label{tab-cost}
\end{table}

\subsection{Convergence result}  
It is easy to see from (\ref{v}) that for any $V\in \hat\partial^2 \phi(w)$, there is
\[
V-I\succeq 0,
\]
implying that $V$ is positive definite. Also note that $\partial^2\phi(w)\subseteq\hat\partial^2 \phi(w)$, we have the following proposition.

\begin{proposition}\label{prop}
	For any $V\in\hat\partial^2 \phi(w)$, $V$ is positive definite.
\end{proposition}

With Proposition \ref{prop}, we are ready to give the local convergence result of semismooth {Newton's} method \ref{algorithm2}.

\begin{theorem}\cite[Thm.3.2]{QiSun93} \label{thm-1} {Let $\omega^*$
		be  a solution of (\ref{sub}). Then every sequence generated by  (\ref{algorithm2})  is superlinearly convergent to $\omega^*$, provided that the starting point $\omega^0$ is sufficiently close to $\omega^*$. Moreover, 
		the convergence rate is quadratic.}
\end{theorem}
\begin{remark}
Recall that in ALM, when solving subproblem, the convergence analysis requires the stopping criteria (A), (B1) and (B2) to be satisfied.
In practice, as  mentioned in \cite[Page 446-447]{LiSun2017}, when $\nabla\phi_k(w^{k+1})$ is sufficiently small, the stopping criteria (A), (B1) and (B2) will be satisfied.
\end{remark}

\section{Algorithm for $\epsilon$-L1-Loss SVR}\label{sec-svr}
In this part, we briefly discuss applying ALM developed in Section \ref{sec3} and \ref{sec4} to $\epsilon$-L1-loss SVR (\ref{prob-svr}).

By letting
$$
B=\left[
\begin{array}{c}
x_1^T\\
\vdots\\
x_m^T
\end{array}
\right],\ y=\left[
\begin{array}{c}
y_1\\
\vdots\\
y_m
\end{array}
\right],
$$
we reformulate (\ref{prob-svr}) as
\be\label{svr-1}
\begin{aligned}
	&\min_w\frac{1}{2}\|w\|^2+{C\sum_{i =1}^m\max(0,|s_i|-\epsilon)} \\
	&{\rm s.t.}\ s=Bw-y.
\end{aligned}
\ee
ALM then can be applied to solve (\ref{svr-1}) with
$$
L_{\sigma}(w,s;\lambda)=\frac{1}{2}\|w\|^2+p_{\epsilon}(s)-\langle\lambda,s-(Bw-y)\rangle+\frac{\sigma}{2}\|s-Bw+y\|^2$$
where $p_\epsilon(\cdot)$ is defined as in (\ref{p-e}). When applying semismooth Newton's method, there is 
\[
\varphi(w) =\min_s\frac{\sigma}{2}\|s-z(w)\|^2+p_{\epsilon}(s)+\frac{1}{2}\|w\|^2-\frac{1}{2\sigma}\|\lambda\|^2
\]
with the gradient
\[
\nabla\varphi(w) = w+B^T\lambda-\sigma B^T({\rm Prox}^{{1/\sigma}} _{p_{\epsilon}}(z{(w)})-Bw+y)
\]
and
\[
\partial^2\varphi(w)(h)\subseteq\hat\partial^2\varphi(\omega)(h),\quad \forall\ h\in \Rm^n,
\]
where $\hat\partial^2\varphi(\omega)=I+\sigma B^TB-\sigma B^T\partial {\rm Prox} ^{1/\sigma}_{p_{\epsilon}}(z(w))B$
and
${z(w)=Bw-y+\frac{\lambda}{\sigma}}$. Recall that $\partial {\rm Prox} ^{1/\sigma}_{p_{\epsilon}}(z(w))$ is the Clarke subdifferential of $\hbox{Prox} ^{1/\sigma}_{p_{\epsilon}}(\cdot)$ at $z(w)$.

When applying CG to solve the linear system, we can choose $V\in\hat\partial^2\varphi(\omega)$ in the following way
\[{
	Vh=h+\sigma B^T(I-U)Bh
	=h+\sigma B(\mathcal{I}(z),:)^TB(\mathcal{I}(z),:)h, \quad U\in\partial  {\rm Prox} ^{1/\sigma}_{p_{\epsilon}}(z(w))},
\]
where \[{\mathcal{I}}(z)=\{i:z_i\in (\epsilon,\epsilon+\frac{C}{\sigma})\cup(-\epsilon-\frac{C}{\sigma},-\epsilon)\}.
\]

\begin{remark}
	For ALM solving $\epsilon$-L1-Loss SVR, the convergence results in Theorems \ref{thm-global} and  \ref{thm-local} (i) holds for ALM. It is not clear whether Theorem \ref{thm-local} (ii) holds or not. The reason is as follows. As shown in Proposition \ref{p_conjugate_svr}, $p_{\epsilon}^*(y)$ equals to $\epsilon\|y\|_1$ if $y\in \overline{P}$ and $+\infty$ otherwise. Hence $p_{\epsilon}(y)$ is neither an indicator function $\delta_M(\cdot)$ or a support function $\delta_M^*(\cdot)$ for some non-empty polyhedral convex set $M\subset\Rm^m$. Therefore, assumption 2.5 in \cite{LiSun2017} fails, putting the metric subregularity of $T_l$ a question since Theorem 2.7 in \cite{LiSun2017} may not hold. As a result, Theorem \ref{thm-local} (ii) may not hold for ALM when solving $\epsilon$-L1-Loss SVR.
\end{remark}

\section{Numerical results}\label{sec5}

In this section, we will conduct numerical tests on various data collected from LIBLINEAR to show the performance of ALM. It is divided into  four parts.  In the first part, we will test the performance of  ALM from various aspects. In the second part, we will test our method based on different parameters. In the third part, we will compare  with one of the competitive solvers DCD (Dual Coordinate Descent method) \cite{LinLiblinear}  in LIBLINEAR for L1-loss SVC. In the final part, we  test our algorithm for $\epsilon$-L1-loss SVR {in comparison with DCD \cite{DCDSVR} in LIBLINEAR}.

{\bf Implementations.} Our algorithm is denoted as ALM-SNCG (Augmented Lagrangian method with semismooth Newton-CG as subsolver). To speed up the solver, we use the following calculation for $\hbox{Prox}^M_p(z)$ in (\ref{Prox-svc}) and $\hbox{Prox}^M_{p_\epsilon}(z)$ respectively,
\[
\begin{array}{ll}
\hbox{Prox}_p^M(z) &= \max(z-CM,0) +
\min(z,0), \\ \hbox{Prox}^M_{p_\epsilon}(z)
&= \max(\min(z,\max(z-CM,\epsilon)), \min(z+CM,-\epsilon)).
\end{array}
\]

All the numerical tests are conducted in Matlab R2017b in Windows 10 on a  {Dell Laptop with an Intel(R) Core(TM) i7-5500U CPU at 2.40GHz and 8 GB of RAM}.  All the data are collected from LIBLINEAR which can be downloaded from {https://www.csie.ntu.edu.tw/~cjlin/libsvmtools/datasets}. The information of datasets in LIBLINEAR is summarized in Table \ref{tab-1} for SVC and Table \ref{tab-svr-info} for SVR.

We will report the following information: the number of iterations in ALM $k$, {the number of iterations in semismooth Newton's method $j$,} the total number of iterations for semismooth Newton's method $it_{sn}$, the total number of iterations in CG $it_{cg}$, the {cpu-time}  $t$ in second. For SVC, we report  accuracy for prediction, which is calculated by
\[\frac{number\ of\ correct\ prediction}{number\ of\ test\ data}\times100\%.\]
For SVR, we report mean squared error (mse) for prediction, which is given by
{
	\[
	mse=\frac{1}{l}\sum_{i=1}^{l}(y_i-\hat{y}_i)^2,
	\]
	where $\hat{y}_i=w^Tx_i$ is the observed value corresponding to the testing data $x_i$, $i=1,\cdots,l$.
}

\begin{table}[H]
	\centering
	\normalsize
	\tbl{Data Information. $m$ is the number of samples, $n$ is the number of features, {``nonzeros" indicates} the number of non-zero elements in all training data, and {``density"} represents the ratio: nonzeros/(m$\times$n)}
	{\begin{tabular}{cccc}
			\hline
			dataset        & ($m$,$n$)           & nonzeros & density  \\ \hline
			leukemia       & (38,7129)       & 270902   & 100.00\% \\
			a1a            & (30956,123)     & 429343   & 11.28\%  \\
			a2a            & (30296,123)     & 420188   & 11.28\%  \\
			a3a            & (29376,123)     & 407430   & 11.28\%  \\
			a4a            & (27780,123)     & 385302   & 11.28\%  \\
			a5a            & (26147,123)     & 362653   & 11.28\%  \\
			a6a            & (21341,123)     & 295984   & 11.28\%  \\
			a7a            & (16461,123)     & 228288   & 11.28\%  \\
			a8a            & (22696,123)     & 314815   & 11.28\%  \\
			a9a            & (32561,123)     & 451592   & 11.28\%  \\
			w1a            & (47272,300)     & 551176   & 3.89\%   \\
			w2a            & (46279,300)     & 539213   & 3.88\%   \\
			w3a            & (44837,300)     & 522338   & 3.88\%   \\
			w4a            & (42383,300)     & 493583   & 3.88\%   \\
			w5a            & (39861,300)     & 464466   & 3.88\%   \\
			w6a            & (32561,300)     & 379116   & 3.88\%   \\
			w7a            & (25057,300)     & 291438   & 3.88\%   \\
			w8a            & (49749,300)     & 579586   & 3.88\%   \\
			breast-cancer  & (683,10)        & 6830     & 100.00\% \\
			cod-rna        & (59535,8)       & 476273   & 100.00\% \\
			diabetes       & (768,8)         & 6135     & 99.85\%  \\
			fourclass      & (862,2)         & 1717     & 99.59\%  \\
			german.numer   & (1000,24)       & 23001    & 95.84\%  \\
			heart          & (270,13)        & 3378     & 96.24\%  \\
			australian     & (690,14)        & 8447     & 87.44\%  \\
			ionosphere     & (351,34)        & 10551    & 88.41\%  \\
			covtype.binary & (581012,54)     & 31363100 & 99.96\%  \\
			ijcnn1         & (49990,22)      & 649869   & 59.09\%  \\
			sonar          & (208,60)        & 12479    & 99.99\%  \\
			splice         & (1000,60)       & 60000    & 100.00\% \\
			svmguide1      & (3089,4)        & 12356    & 100.00\% \\
			svmguide3      & (1243,22)       & 22014    & 80.50\%  \\
			phishing       & (11055,68)      & 751740   & 100.00\% \\
			madelon        & (2000,500)      & 979374   & 97.94\%  \\
			mushrooms      & (8124,112)      & 901764   & 99.11\%  \\
			duke breast-cancer          & (44,7129)       & 313676   & 100.00\% \\
			gisette        & (6000,5000)     & 29729997 & 99.10\%  \\
			news20.binary  & (19996,1355191) & 9097916  & 0.03\%   \\
			rcv1.binary    & (20242,47236)   & 1498952  & 0.16\%   \\
			real-sim       & (72309,20958)   & 3709083  & 0.24\%   \\
			livery         & (145,5)         & 724      & 99.86\%  \\
			colon-cancer   & (62,2000)       & 124000   & 100.00\% \\
			skin-nonskin   & (245057,3)      & 735171   & 100.00\% \\ \hline
	\end{tabular}}
	\label{tab-1}
\end{table}

\begin{table}[H]
	\centering
	\tbl{{Data Information for SVR}}
	{\begin{tabular}{cccc}
			\hline
			dataset     & ($m$,$n$)          & nonzeros & density  \\ \hline
			abalone     & (4177,8)       & 32080    & 96.00\%  \\
			bodyfat     & (252,14)       & 3528     & 100.00\% \\
			cadata      & (20640,8)      & 165103   & 99.99\%  \\
			cpusmall    & (8192,12)      & 98304    & 100.00\% \\
			E2006-train & (16087,150360) & 19971014 & 0.83\%   \\
			E2006-test  & (3308,150358)  & 4559527  & 0.92\%   \\
			eunite2001  & (336,16)       & 2651     & 49.31\%  \\
			housing     & (506,13)       & 6578     & 100.00\% \\
			mg          & (1385,6)       & 8310     & 100.00\% \\
			mpg         & (392,7)        & 2641     & 96.25\%  \\
			pyrim       & (74,27)        & 1720     & 86.09\%  \\
			space-ga    & (3107,6)       & 18642    & 100.00\% \\
			triazines   & (186,60)       & 9982     & 89.44\%  \\ \hline
	\end{tabular}}
	\label{tab-svr-info}
\end{table}

\subsection{Performance test}
In this part, we will show the performance of ALM-SNCG  including the low computational cost, the convergence rate of semismooth Newton's method for solving subproblems as well as the accuracy with respect to the iterations in ALM.

\begin{figure}[H]
	\centering
	\includegraphics[width=15cm]{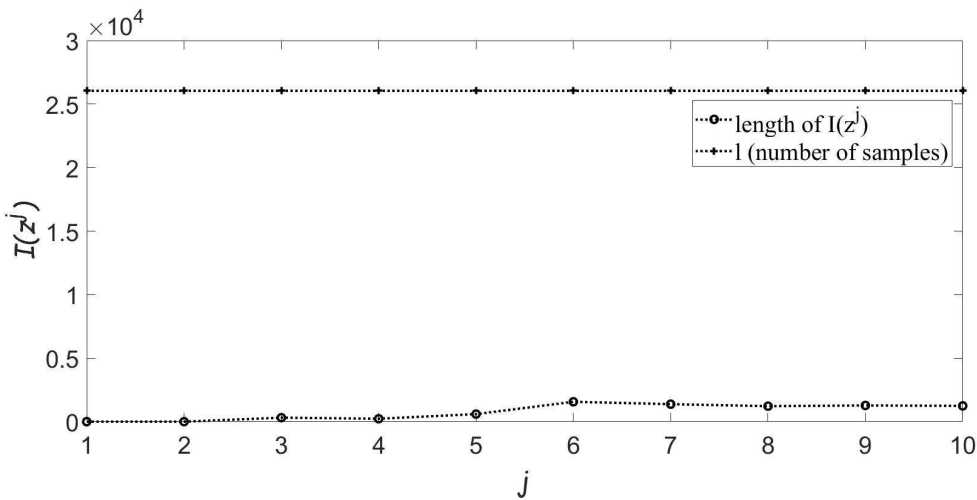}
	\caption{a9a: $|\mathcal{I}(z^j)|$}
	\label{fig-1}
\end{figure}

{\bf Low computational cost.} As we show in Section \ref{sec4-2}, we can reduce the computational cost significantly by exploring the sparse structure when calculating $Vh$ in CG. {Without causing any chaos, we denote $z(w^j)$ as $z^j$.} Below we plot {$|\mathcal{I}(z^j)|$} with respect to $j$ for data {a9a} ($m = 32,561, n = 123$) in the first loop of ALM. We have \[|\mathcal{I}(z^j)|\ =[ 0,\ 0,\ 314,\ 229,\ 595,\ 1564,\ 1377,\ 1219,\ 1275,\ 1241],\]
and the radio of $|\mathcal{I}(z^j)|$ over $m$ is \[[0,\ 0,\ 0.0096,\ 0.0070,\ 0.0183,\ 0.0480,\ 0.0423,\ 0.0374,\ 0.0392,\ 0.0381].\]
From Figure $\ref{fig-1}$, one can see that compared with the large sample size $m = 32,561$, the number of elements in $\mathcal{I}(z)$ is far more less than $m$, which only accounts for {less than $5\%$} of $m$. Consequently, this verifies our claim that the computational complexity in WAY III is significantly reduced.


\begin{figure}[H]
	\centering
	\includegraphics[width=15cm]{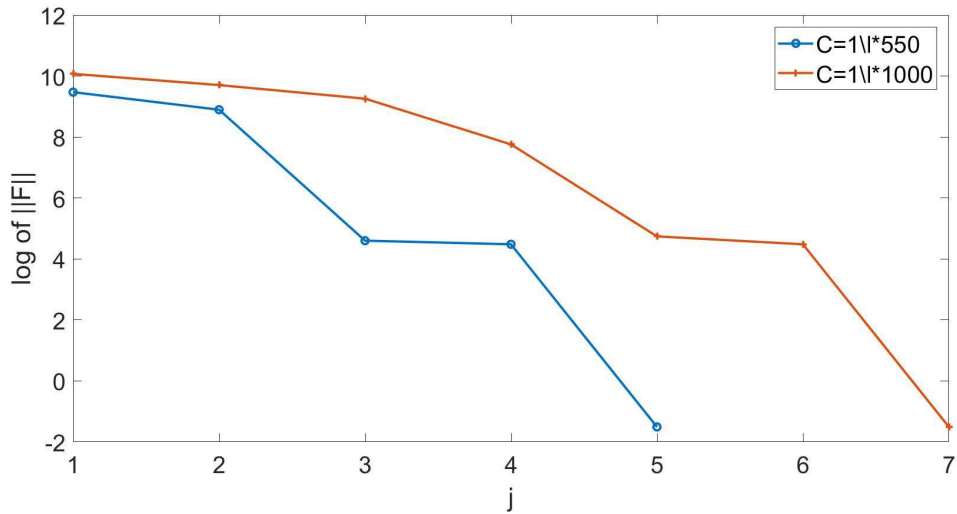}
	\caption{$\log\|F(w^j)\|$  with C=$1/l*550$ and $1/l*1000$}
	\label{fig-2}
\end{figure}

{\bf Convergence of semismooth Newton's method.} To demonstrate the quadratic convergence of semismooth Newton's method, we plot $\|F(w^j)\|$ with respect to $j$ in the first loop of ALM for dataset leukemia. The result is demonstrated in Figure \ref{fig-2}, where the quadratic convergence rate can be indeed observed.


\begin{figure}[H]
	\centering
	\includegraphics[width=15cm]{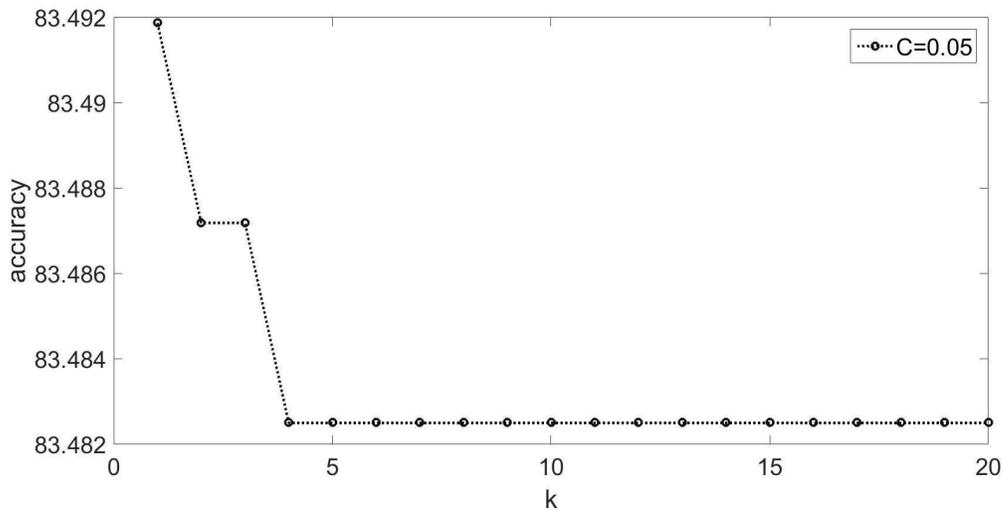}
	\caption{Accuracy during the update of ALM for a6a with $ C=0.05$.}
	\label{fig-3}
\end{figure}

{\bf Accuracy.} Note that we are solving a practical problem, with predicting purpose. We are interested in the following question: as the iteration keeps on in ALM, what it happens to the accuracy? Below, we plot the accuracy vs outer iteration number $k$ in Figure \ref{fig-3}. As we can see,  the accuracy is reduced significantly during the first few iterations in ALM. However, little progress has been made after that. Consequently, we set the maximum number of outer iterations  as $10$ in our test.


\subsection{Performance with different parameters}
It is noted that the choice of $C$ in L1-loss SVC  (\ref{prob}) may have great influence on the performance of the algorithms. Consequently, in this part, we test our algorithms with different $C$s. To speed up the code, we  try {two choices}. {One is to} generate starting point by solving subproblem {(\ref{prob-sub})} using {alternating direction method of multipliers} (ADMM). {The other is to use a fixed starting point}  {$w^0=ones(n,1)$}.  We report the results in Table \ref{tab-2} and \ref{tab-3}.
One can see that {ALM-SNCG} {with fixed starting point} has a better accuracy  than that with ADMM as starting point. For {cpu-time},  using ADMM or not does not make much difference. 
For the choice of C, it seems that both values of $C=1/l*550,\ 1/l*1000$ give comparable accuracy and {cpu-time}.  In our following test, we choose $C=1/l*550$ and {use $w^0=ones(n,1)$ as our starting point}.

\subsection{Comparison with LIBLINEAR for SVC}
As pointed out in \cite{review18}, LIBLINEAR is the most popular solver. Consequently, in this part, we compare our algorithm with solvers in LIBLINEAR on dataset for SVC. We split each set to 80\% training and 20\% testing. Specifically, we compare with DCD which solves also model (\ref{prob}).
Parameters in ALM are set as $\tau=1$, $\sigma_0=0.15$, $\sigma_{\max}=2$, $\theta=0.8$ and $\lambda_0=zeros(m,1)$.
For semismooth Newton's method, we choose
$\rho=0.5$,  $\delta=\max(1.0\times10^{-2},\ 1.0\times10^{-j})$. 
The results are reported in Table \ref{tab-4}.
%
%



\begin{table}[H]
	\centering
	\tbl{Comparison of Starting Points: C=1/l*550. I: $w^0=ones(n,1)$ as starting point; II: with ADMM as starting point}
	{\begin{tabular}{ccccccc}
		\hline
		\multirow{2}{*}{dataset} & \multicolumn{2}{c}{iter($k,it_{sn},it_{cg}$)} & \multicolumn{2}{c}{times(s)} & \multicolumn{2}{c}{accuracy($\%$)}  \\
		& I           & II                 & I       & II         & I          & II             \\ \hline
		leukemia                            & (1,5,11)             & (1,5,11)             & 0.015         & 0.097        & \textbf{100.000} & \textbf{100.000} \\
		a1a                                 & (1,11,77)            & (1,11,58)            & 0.036         & 0.033        & \textbf{84.819}  & 84.787           \\
		a2a                                 & (1,11,80)            & (1,11,64)            & 0.029         & 0.031        & \textbf{84.818}  & 84.769           \\
		a3a                                 & (1,9,70)             & (1,11,63)            & 0.024         & 0.028        & \textbf{84.717}  & \textbf{84.717}  \\
		a4a                                 & (1,10,73)            & (1,11,64)            & 0.028         & 0.027        & 84.737           & \textbf{84.773}  \\
		a5a                                 & (1,10,78)            & (1,12,70)            & 0.022         & 0.034        & \textbf{84.780}  & 84.761           \\
		a6a                                 & (1,9,67)             & (1,10,51)            & 0.022         & 0.021        & \textbf{84.587}  & \textbf{84.587}  \\
		a7a                                 & (1,11,89)            & (1,10,70)            & 0.019         & 0.018        & \textbf{84.604}  & 84.573           \\
		a8a                                 & (1,9,72)             & (1,10,55)            & 0.019         & 0.020        & \textbf{84.031}  & 83.965           \\
		a9a                                 & (1,11,84)            & (1,11,66)            & 0.031         & 0.030        & 84.677           & \textbf{84.692}  \\
		w1a                                 & (1,14,99)            & (1,17,114)           & 0.057         & 0.072        & \textbf{99.958}  & \textbf{99.958}  \\
		w2a                                 & (1,12,84)            & (1,14,122)           & 0.062         & 0.061        & \textbf{99.968}  & 99.957           \\
		w3a                                 & (1,17,141)           & (1,20,195)           & 0.093         & 0.148        & \textbf{99.955}  & \textbf{99.955}  \\
		w4a                                 & (1,13,95)            & (1,12,91)            & 0.059         & 0.067        & \textbf{99.953}  & \textbf{99.953}  \\
		w5a                                 & (2,22,168)           & (1,16,138)           & 0.105         & 0.087        & \textbf{99.950}  & \textbf{99.950}  \\
		w6a                                 & (1,17,170)           & (1,20,151)           & 0.068         & 0.065        & \textbf{99.923}  & \textbf{99.923}  \\
		w7a                                 & (1,14,100)           & (1,12,90)            & 0.031         & 0.033        & \textbf{99.900}  & \textbf{99.900}  \\
		w8a                                 & (1,14,114)           & (1,14,125)           & 0.060         & 0.062        & \textbf{99.950}  & \textbf{99.950}  \\
		breast-cancer                       & (11,32,116)          & (11,32,108)          & 0.010         & 0.005        & \textbf{99.270}  & \textbf{99.270}  \\
		cod-rna                             & (1,10,32)            & (1,9,31)             & 0.032         & 0.036        & \textbf{20.417}  & 20.324           \\
		diabetes                            & (11,35,114)          & (11,33,107)          & 0.006         & 0.004        & \textbf{75.974}  & \textbf{75.974}  \\
		fourclass                           & (11,30,50)           & (11,28,47)           & 0.004         & 0.002        & \textbf{76.879}  & \textbf{76.879}  \\
		german.numer                        & (11,44,305)          & (11,39,272)          & 0.020         & 0.016        & \textbf{78.500}  & \textbf{78.500}  \\
		heart                               & (11,27,132)          & (11,28,141)          & 0.006         & 0.005        & \textbf{85.185}  & \textbf{85.185}  \\
		australian                          & (11,37,232)          & (11,36,216)          & 0.009         & 0.007        & \textbf{86.232}  & \textbf{86.232}  \\
		ionosphere                          & (11,33,189)          & (11,33,191)          & 0.006         & 0.010        & \textbf{98.592}  & \textbf{98.592}  \\
		covtype.binary                      & (2,68,379)           & (1,46,377)           & 7.956         & 7.883        & \textbf{64.331}  & 64.210           \\
		ijcnn1                              & (1,12,80)            & (1,9,36)             & 0.091         & 0.044        & \textbf{90.198}  & \textbf{90.198}  \\
		sonar                               & (11,31,171)          & (11,31,171)          & 0.009         & 0.009        & \textbf{64.286}  & \textbf{64.286}  \\
		splice                              & (11,42,230)          & (11,41,222)          & 0.021         & 0.023        & \textbf{80.500}  & \textbf{80.500}  \\
		svmguide1                           & (11,30,70)           & (11,30,79)           & 0.007         & 0.006        & \textbf{78.479}  & \textbf{78.479}  \\
		svmguide3                           & (11,34,134)          & (11,30,120)          & 0.013         & 0.013        & \textbf{0.000}   & \textbf{0.000}   \\
		phishing                            & (1,14,101)           & (1,12,76)            & 0.047         & 0.043        & 92.492           & \textbf{92.583}  \\
		madelon                             & (11,51,445)          & (11,46,453)          & 0.743         & 0.864        & \textbf{52.250}  & \textbf{52.250}  \\
		mushrooms                           & (2,34,393)           & (1,25,361)           & 0.202         & 0.220        & \textbf{100.000} & \textbf{100.000} \\
		duke breast-cancer                               & (1,7,14)             & (1,6,13)             & 0.022         & 0.132        & \textbf{100.000} & \textbf{100.000} \\
		gisette                             & (4,187,1174)         & (3,117,1570)         & 20.870        & 23.454       & \textbf{97.583}  & 97.500           \\
		news20.binary                       & (11,29,61)           & (11,29,61)           & 4.463         & 4.389        & \textbf{11.275}  & \textbf{11.275}  \\
		rcv1-train.binary                   & (11,61,142)          & (11,36,111)          & 0.537         & 0.479        & \textbf{94.394}  & \textbf{94.394}  \\
		real-sim                            & (1,32,43)            & (1,12,33)            & 0.436         & 0.406        & \textbf{80.729}  & 80.556           \\
		liver-disorders                     & (11,22,44)           & (11,22,42)           & 0.005         & 0.002        & \textbf{55.172}  & \textbf{55.172}  \\
		colon-cancer                        & (1,7,16)             & (1,7,15)             & 0.009         & 0.014        & \textbf{69.231}  & \textbf{69.231}  \\
		skin-nonskin                        & (1,8,20)             & (1,9,18)             & 0.157         & 0.170        & \textbf{90.429}  & 90.413           \\ \hline
	\end{tabular}}
\label{tab-2}
\end{table}

\begin{table}[H]
	\centering
	\tbl{Comparison of Starting Points: C=1/l*1000}
	{\begin{tabular}{ccccccc}
			\hline
			\multirow{2}{*}{dataset} & \multicolumn{2}{c}{iter($k,it_{sn},it_{cg}$)} & \multicolumn{2}{c}{time(s)}      & \multicolumn{2}{c}{accuracy($\%$)}  \\
			& I           & II                 & I       & II         & I          & II               \\ \hline
			leukemia                        & (1,7,13)             & (1,6,12)             & \textbf{0.027} & 0.121           & \textbf{100.000} & \textbf{100.000} \\
			a1a                        & (8,36,284)           & (8,31,227)           & 0.099          & \textbf{0.079}  & \textbf{84.868}  & \textbf{84.868}  \\
			a2a                        & (10,33,302)          & (10,40,280)          & 0.107          & \textbf{0.096}  & \textbf{84.835}  & \textbf{84.835}  \\
			a3a                        & (11,41,346)          & (11,39,308)          & 0.098          & \textbf{0.095}  & \textbf{84.717}  & \textbf{84.717}  \\
			a4a                        & (11,35,304)          & (11,43,305)          & \textbf{0.100} & 0.107           & \textbf{84.791}  & \textbf{84.791}  \\
			a5a                        & (10,33,280)          & (10,35,263)          & 0.080          & \textbf{0.079}  & \textbf{84.665}  & \textbf{84.665}  \\
			a6a                        & (10,33,294)          & (11,42,291)          & \textbf{0.080} & 0.093           & \textbf{84.657}  & \textbf{84.657}  \\
			a7a                        & (11,36,302)          & (10,37,275)          & \textbf{0.065} & 0.066           & \textbf{84.604}  & \textbf{84.604}  \\
			a8a                        & (11,42,320)          & (11,37,295)          & \textbf{0.075} & 0.077           & \textbf{84.141}  & \textbf{84.141}  \\
			a9a                        & (9,38,322)           & (9,38,270)           & 0.109          & \textbf{0.096}  & \textbf{84.738}  & \textbf{84.738}  \\
			w1a                        & (1,18,157)           & (2,18,135)           & 0.075          & \textbf{0.072}  & 99.905           & \textbf{99.926}  \\
			w2a                        & (1,13,97)            & (1,16,157)           & \textbf{0.054} & 0.080           & \textbf{99.924}  & \textbf{99.924}  \\
			w3a                        & (1,15,142)           & (1,17,123)           & \textbf{0.059} & 0.061           & \textbf{99.944}  & 99.933           \\
			w4a                        & (1,15,132)           & (1,20,156)           & \textbf{0.082} & 0.108           & \textbf{99.929}  & \textbf{99.929}  \\
			w5a                        & (1,14,102)           & (1,13,89)            & 0.059          & \textbf{0.048}  & 99.887           & \textbf{99.900}  \\
			w6a                        & (1,19,148)           & (1,15,109)           & 0.058          & \textbf{0.049}  & \textbf{99.893}  & 99.877           \\
			w7a                        & (1,17,131)           & (5,27,219)           & \textbf{0.047} & 0.066           & 99.900           & \textbf{99.920}  \\
			w8a                        & (1,19,171)           & (1,21,157)           & 0.104          & \textbf{0.097}  & \textbf{99.920}  & \textbf{99.920}  \\
			breast-cancer              & (11,33,130)          & (11,32,124)          & \textbf{0.004} & 0.005           & \textbf{99.270}  & \textbf{99.270}  \\
			cod-rna                    & (1,17,44)            & (1,11,35)            & 0.068          & \textbf{0.045}  & \textbf{21.542}  & 21.492           \\
			diabetes                   & (11,33,102)          & (11,31,96)           & 0.005          & \textbf{0.004}  & \textbf{74.675}  & \textbf{74.675}  \\
			fourclass                  & (11,30,49)           & (11,28,46)           & \textbf{0.002} & 0.005           & \textbf{76.879}  & \textbf{76.879}  \\
			german.numer               & (11,42,305)          & (11,39,304)          & \textbf{0.017} & 0.019           & \textbf{77.500}  & \textbf{77.500}  \\
			heart                      & (11,24,118)          & (11,25,127)          & 0.004          & \textbf{0.003}  & \textbf{88.889}  & \textbf{88.889}  \\
			australian                 & (11,34,195)          & (11,32,171)          & \textbf{0.006} & \textbf{0.006}  & \textbf{86.957}  & \textbf{86.957}  \\
			ionosphere                 & (11,30,158)          & (11,29,151)          & \textbf{0.006} & 0.007           & \textbf{97.183}  & \textbf{97.183}  \\
			covtype.binary             & (3,120,489)          & (3,107,505)          & 14.295         & \textbf{13.708} & 64.738           & \textbf{64.769}  \\
			ijcnn1                     & (11,34,221)          & (11,33,205)          & 0.267          & \textbf{0.208}  & \textbf{90.378}  & \textbf{90.378}  \\
			sonar                      & (11,30,172)          & (11,30,172)          & \textbf{0.022} & 0.057           & \textbf{64.286}  & \textbf{64.286}  \\
			splice                     & (11,41,212)          & (11,39,198)          & 0.031          & \textbf{0.028}  & \textbf{81.500}  & \textbf{81.500}  \\
			svmguide1                  & (11,32,91)           & (11,32,96)           & 0.011          & \textbf{0.008}  & \textbf{86.408}  & \textbf{86.408}  \\
			svmguide3                  & (11,32,149)          & (11,31,152)          & 0.023          & \textbf{0.014}  & \textbf{0.000}   & \textbf{0.000}   \\
			phishing                   & (10,47,363)          & (10,41,305)          & 0.205          & \textbf{0.117}  & \textbf{92.492}  & \textbf{92.492}  \\
			madelon                    & (11,67,478)          & (11,44,420)          & 1.148          & \textbf{1.068}  & \textbf{50.000}  & \textbf{50.000}  \\
			mushrooms                  & (1,31,244)           & (1,26,258)           & \textbf{0.125} & 0.171           & \textbf{100.000} & \textbf{100.000} \\
			duke breast-cancer                      & (1,8,15)             & (1,9,16)             & \textbf{0.024} & 0.103           & \textbf{100.000} & \textbf{100.000} \\
			gisette                    & (6,266,1270)         & (5,235,1159)         & 26.389         & \textbf{20.597} & 97.250           & \textbf{97.500}  \\
			news20.binary              & (11,32,78)           & (11,32,78)           & 5.769          & \textbf{5.760}  & \textbf{49.850}  & \textbf{49.850}  \\
			rcv1-train.binary          & (10,111,194)         & (11,35,125)          & 1.302          & \textbf{0.613}  & \textbf{95.357}  & \textbf{95.357}  \\
			real-sim                   & (11,72,119)          & (1,16,49)            & 1.050          & \textbf{1.007}  & \textbf{74.582}  & 74.063           \\
			liver-disorders            & (11,22,44)           & (11,22,42)           & 0.011          & \textbf{0.005}  & \textbf{55.172}  & \textbf{55.172}  \\
			colon-cancer               & (1,6,14)             & (1,6,18)             & \textbf{0.007} & 0.011           & \textbf{69.231}  & \textbf{69.231}  \\
			skin-nonskin               & (1,11,27)            & (1,13,24)            & \textbf{0.189} & 0.218           & \textbf{90.235}  & 90.231           \\ \hline
	\end{tabular}}
	\label{tab-3}
\end{table}

\begin{table}[H]
	\centering
	\normalsize
	\tbl{The comparison results for L1-loss SVC. }
	{\begin{tabular}{ccc}
		\hline
		dataset           & {time(s)(DCD$|$ALM-SNCG)} & {accuracy($\%$)(DCD$|$ALM-SNCG)}               \\ \hline
		leukemia               & 0.017 $|$ \textbf{0.015}   & \textbf{100.000} $|$ \textbf{100.000}                       \\
		a1a               & 0.059 $|$ \textbf{0.036}   & 84.722 $|$ \textbf{84.819}                        \\
		a2a               & 0.049 $|$ \textbf{0.029}   & 84.752 $|$ \textbf{84.818}                        \\
		a3a               & 0.044 $|$ \textbf{0.024}   & 84.700 $|$ \textbf{84.717}                        \\
		a4a               & 0.042 $|$ \textbf{0.028}   & 84.683 $|$ \textbf{84.737}                        \\
		a5a               & 0.040 $|$ \textbf{0.022}   & 84.627 $|$ \textbf{84.780}                        \\
		a6a               & 0.033 $|$ \textbf{0.022}   & 84.516 $|$ \textbf{84.587}                        \\
		a7a               & 0.026 $|$ \textbf{0.019}   & 84.482 $|$ \textbf{84.604}                        \\
		a8a               & 0.033 $|$ \textbf{0.019}   & \textbf{84.119} $|$ 84.031                                 \\
		a9a               & 0.068 $|$ \textbf{0.031}   & \textbf{84.708} $|$ 84.677                                 \\
		w1a               & 0.088 $|$ \textbf{0.057}   & 99.662 $|$ \textbf{99.958}                        \\
		w2a               & 0.071 $|$ \textbf{0.062}   & 99.676 $|$ \textbf{99.968}                        \\
		w3a               & 0.103 $|$ \textbf{0.093}   & 99.699 $|$ \textbf{99.955}                        \\
		w4a               & 0.077 $|$ \textbf{0.059}   & 99.705 $|$ \textbf{99.953}                        \\
		w5a               & 0.121 $|$ \textbf{0.105}   & 99.737 $|$ \textbf{99.950}                        \\
		w6a               & \textbf{0.047} $|$ 0.068   & 99.754 $|$ \textbf{99.923}                        \\
		w7a               & 0.039 $|$ \textbf{0.031}   & 99.721 $|$ \textbf{99.900}                        \\
		w8a               & 0.094 $|$ \textbf{0.060}   & 99.668 $|$ \textbf{99.950}                        \\
		breast-cancer     & \textbf{0.001} $|$ 0.010   & \textbf{100.000} $|$ 99.270                                 \\
		cod-rna           & 0.066 $|$ \textbf{0.032}   & 19.879 $|$ \textbf{20.417}                        \\
		diabetes          & \textbf{0.001} $|$ 0.006            & 74.675 $|$ \textbf{75.974}                        \\
		fourclass         & \textbf{0.001} $|$ 0.004            & 68.786 $|$ {\color[HTML]{FE0000} \textbf{76.879}} \\
		german.numer      & \textbf{0.004} $|$ 0.020  & 78.000 $|$ \textbf{78.500}                        \\
		heart             & \textbf{0.002} $|$ 0.006            & 81.481 $|$ \textbf{85.185}                        \\
		australian        & \textbf{0.001} $|$ 0.009            & \textbf{86.232} $|$ \textbf{86.232}                        \\
		ionosphere        & \textbf{0.006} $|$ \textbf{0.006}   & \textbf{98.592} $|$ \textbf{98.592}                        \\
		covtype.binary    & 12.086 $|$ \textbf{7.956}   & 64.094 $|$ \textbf{64.331}                        \\
		ijcnn1            & 0.125 $|$ \textbf{0.091}   & \textbf{90.198} $|$ \textbf{90.198}                        \\
		sonar             & \textbf{0.006} $|$ 0.009            & \textbf{64.286} $|$ \textbf{64.286}                        \\
		splice            & \textbf{0.016} $|$ 0.021            & 69.500 $|$ {\color[HTML]{FE0000} \textbf{80.500}} \\
		svmguide1         & \textbf{0.001} $|$ 0.007            & 53.560 $|$ {\color[HTML]{FE0000} \textbf{78.479}} \\
		svmguide3         & \textbf{0.002} $|$ 0.013            & \textbf{0.000} $|$ \textbf{0.000}                         \\
		phishing          & 0.081 $|$ \textbf{0.047}   & \textbf{92.763} $|$ 92.492                                 \\
		madelon           & \textbf{0.181} $|$ 0.743            & \textbf{52.500} $|$ 52.250                                 \\
		mushrooms         & \textbf{0.131} $|$ 0.202            & \textbf{100.000} $|$ \textbf{100.000}                       \\
		duke breast-cancer             & 0.034 $|$ \textbf{0.022}   & \textbf{100.000} $|$ \textbf{100.000}                       \\
		gisette           & \textbf{3.796} $|$ 19.270           & 97.333 $|$ \textbf{97.583}                        \\
		news20.binary     & \textbf{0.642} $|$ 4.463            & \textbf{27.775} $|$ 11.275                                 \\
		rcv1-train.binary & \textbf{0.123} $|$ 0.537            & \textbf{94.443} $|$ 94.394                                 \\
		real-sim          & \textbf{0.353} $|$ 0.436            & 72.383 $|$ {\color[HTML]{FE0000} \textbf{80.729}} \\
		liver-disorders   & \textbf{0.000} $|$ 0.005            & 41.379 $|$ {\color[HTML]{FE0000} \textbf{55.172}} \\
		colon-cancer      & \textbf{0.006} $|$ 0.009            & \textbf{69.231} $|$ \textbf{69.231}                        \\
		skin-nonskin      & 0.171 $|$ \textbf{0.157}   & 88.236 $|$ \textbf{90.429}                        \\ \hline
	\end{tabular}}
	\label{tab-4}
\end{table}

{We can get the following observations from the results.}
\bit\item  Both of the two algorithms can obtain high accuracy. We marked the winners of {cpu-time} and accuracy in bold. The accuracy of more than 80\% of datasets is over {72\%}, and some datasets's accuracy is more than 90\%. In comparison with DCD for the 43 classification datasets, {ALM-SNCG} has equal or higher accuracy for 36 datasets. In particular, for dataset {fourclass, splice, svmguide1, real-sim, and liver-disorders, the accuracy of {ALM-SNCG} has increased by nearly 10\% or even more than 10\% (marked in red).}

\item In terms of cpu-time,  {ALM-SNCG}  {is competitive} with DCD. For example, for the dataset covtype.binary, {ALM-SNCG} saves {cpu-time} by almost 50\%. For all the 43 datasets, {ALM-SNCG} is faster than DCD for 24 datasets.

\eit

\subsection{Comparison with LIBLINEAR for SVR}
In this part, we compare our algorithm with DCD \cite{DCDSVR} in LIBLINEAR on dataset for SVR (\ref{prob-svr}). {We split each dataset to 60\% training and 40\% testing. }
{We choose C = 1/n*5. Parameters in Algorithm \ref{alg-alm-1} are set as
	$\sigma_0=0.1$  and  other parameters are set as the same in SVC.}
The results are reported in Table \ref{tab-svr-result}, from which we have the following observations.
{
	\bit\item  In comparison with DCD for the 13 {regression} datasets, ALM-SNCG has equal or higher accuracy for all the 13 datasets.
	
	\item In terms of cpu-time,  ALM-SNCG  {is competitive with} DCD. For example, for the dataset E2006-train,  ALM-SNCG  saves {cpu-time} by more than 50\%. {For all the 13 datasets, ALM-SNCG has equal or shorter cpu-time for 10 datasets.}
	
	\eit
	
}

\begin{table}[H]
	\centering
	\tbl{The comparison results for SVR.}
	{\begin{tabular}{ccc}
		\hline
		dataset        & {time(s)(DCD$|$ALM-SNCG)} & {mse(DCD$|$ALM-SNCG)}   \\ \hline
		abalone     & \textbf{0.002} $|$ \textbf{0.002}   & 27.47 $|$ \textbf{12.85}     \\
		bodyfat     & \textbf{0.001} $|$ \textbf{0.001}   & 0.02 $|$ \textbf{0.01}      \\
		cadata      & 0.014 $|$ \textbf{0.008}   & 0.17 $|$ \textbf{0.15}      \\
		cpusmall    & 0.004 $|$ \textbf{0.001}   & 2253.55 $|$ \textbf{1862.67}   \\
		E2006-train & 1.517 $|$ \textbf{0.733}   & 0.31 $|$ \textbf{0.29}      \\
		E2006-test  & \textbf{0.316} $|$ 0.349            & 0.23 $|$ \textbf{0.21}      \\
		eunite2001  & \textbf{0.001} $|$ \textbf{0.001}   & 539999.28 $|$ \textbf{532688.85} \\
		housing     & \textbf{0.001} $|$ \textbf{0.001}   & 134.79 $|$ \textbf{93.77}     \\
		mg          & \textbf{0.000} $|$ 0.001            & 0.42 $|$ \textbf{0.02}      \\
		mpg         & \textbf{0.000} $|$ \textbf{0.000}   & 862.10 $|$ \textbf{606.14}    \\
		pyrim       & \textbf{0.000} $|$ 0.001            & 0.02 $|$ \textbf{0.01}      \\
		space-ga    & \textbf{0.001} $|$ \textbf{0.001}   & 0.26 $|$ \textbf{0.03}      \\
		triazines   & \textbf{0.002} $|$ \textbf{0.002}   & \textbf{0.03} $|$ \textbf{0.03}      \\ \hline
	\end{tabular}}
\label{tab-svr-result}
\end{table}

\section{Conclusions}\label{sec6} 
In this paper, we {proposed} a semismooth Newton-CG based on augmented Lagrangian method for solving the L1-loss SVC and $\epsilon$-L1-loss SVC. 
The proposed algorithm {enjoyed} the traditional convergence result while keeping the fast quadratic convergence and low computational complexity in semismooth Newton algorithm as a subsolver. Extensive numerical results on datasets in LIBLINEAR demonstrated the superior performance of the proposed algorithm over the LIBLINEAR in terms of both accuracy and speed.\\
\\

\noindent{\bf Acknowledgements.}\\
We would like to thank the three anonymous reviewers for there valuable comments which helped improve the paper significantly. We would also like to thank Prof. Chao Ding from Chinese Academy of Sciences for providing the important reference of  Prof. Jie Sun's PhD thesis. We would also like to thank Dr. Xudong Li from Fudan University for helpful discussions.\\
\\

\noindent{\bf Notes on contributors.}\\
\emph{Yinqiao Yan}  received his Bachelor of science in statistics degree from Beijing Institute of Technology,
China, in 2019. He is pursuing his PhD degree in Institute of Statistics and Big Data, Renmin University of China, China.\\

\noindent \emph{Qingna Li} received her Bachelor’s degree in information and computing science and Doctor’s degree
in computational mathematics from Hunan University, China, in 2005 and 2010 respectively. Currently,
she is an associate professor in School of Mathematics and Statistics, Beijing Institute of
Technology. Her research interests include continuous optimization and its applications.

\section*{Appendix. The calculation of $\hbox{Prox}_p^M(z)$}

Note that
\begin{eqnarray*}
	&&\min_s\frac1{2M}\|z-s\|^2+C\sum_{i = 1}^m \max(0,s_i)\\
	&\Longleftrightarrow &\min_s \sum_{i = 1}^m(\frac1{2M}(z_i-s_i)^2+C \max(0,s_i))\\
	&\Longleftrightarrow &\min_{s_i} \frac1{2M}(z_i-s_i)^2+C \max(0,s_i), \ i  = 1,\dots, m.\\
\end{eqnarray*}
For each $i$, consider solving the following problem
\[
\min_{s_i}
\frac1{2M}(z_i-s_i)^2+C \max(0,s_i).
\]
This is a problem of minimizing a piecewise quadratic function. The objective function of this problem takes the following form
\[
\left\{
\begin{aligned}
&\frac1{2M}(z_i-s_i)^2,\ & s_i\le0,\\
&\frac1{2M}s_i^2+(C-\frac1{M}z_i)s_i+\frac1{2M}z_i^2,\ & s_i>0.
\end{aligned}
\right.
\]
One can see that the solution is
\[
s_i^*=\left\{
\begin{array}{ll}
z_i-CM,\ & z_i>CM,\\
z_i,\ & z_i<0,\\
0,\ & 0\leq z_i\leq CM.
\end{array}
\right.
\]
Therefore, we get that
\[
(\hbox{Prox}_p^M(z))_i = \left\{
\begin{array}{ll}
z_i-CM,\ & z_i>CM,\\
z_i,\ & z_i<0,\\
0,\ & 0\leq z_i\leq CM.
\end{array}
\right.
\]

\end{document}